\documentclass[a4paper,10pt]{article}
\usepackage{geometry,times}
\usepackage{amsfonts}
\usepackage{amssymb}
\usepackage{amsmath}
\usepackage{tikz}
\usepackage{pgfplots}
\usepackage{tikz-3dplot}
\definecolor{paleyellow}{rgb}{1, 1, 0.9}
\usepackage[backgroundcolor=paleyellow]{todonotes}
\usepackage{amsthm}
\newcommand{\rose}{R\accent 23 u\v zi\v cka}
\newcommand{\bb}{\mathbb}

\theoremstyle{plain}
\newtheorem{satz}{prop}

\newtheorem{lem}[satz]{Lemma}

\newtheorem{prop}[satz]{Proposition}
\theoremstyle{thm}
\newtheorem{theor}{Theorem}
\theoremstyle{note}

\theoremstyle{definition}

\theoremstyle{problem}

\theoremstyle{notadef}
\newtheorem{notation}{Notation}
\usepackage[]{graphicx}
\usepackage{color}
\definecolor{dark-BLGR}{rgb}{0.28,0.28,0.28}

\definecolor{dark-red}{rgb}{0.90,0.14,0.14}

\definecolor{dark-viol}{rgb}{0.73,0.18,0.69}
\definecolor{Cyan}{rgb}{0.31,0.67,0.82}

\definecolor{light-yellow}{rgb}{1,1,0.8}

\definecolor{Blue}{rgb}{0.0353,0.0275,0.4}

\definecolor{GRAY}{rgb}{0.26,0.26,0.26}

\definecolor{GREYY}{rgb}{0.45,0.45,0.45}

\definecolor{vivid-viol}{rgb}{0.3255,0.0353,0.55}

\definecolor{dark-blue}{rgb}{0.05,0.05,0.65}

\definecolor{dark-green}{rgb}{0.03,0.77,0.29}

\definecolor{dark-Green}{rgb}{0.03,0.57,0.09}
\newcommand{\DarkGr}{\textcolor{dark-Green}}
\definecolor{strong-viol}{rgb}{0.2353,0.094,0.349}

\definecolor{BLUE}{rgb}{0.41,0.44,0.93}

\definecolor{RED}{rgb}{0.90,0.18,0.32}

\definecolor{Yel}{rgb}{0.89,0.95,0.19}

\definecolor{White}{rgb}{1,1,1}
\definecolor{black}{rgb}{0,0,0}
\newcommand{\Black}{\textcolor{black}}
\definecolor{Orchid}{rgb}{0.6,0.1607,0.2}

\definecolor{Orange}{rgb}{0.99,0.49,0}

\definecolor{Thistle}{rgb}{0.1151,0.1249,0.9873}

\definecolor{brown}{rgb}{0.298,0.153,0.07843}

 \usepackage{titlesec}
\titleformat{\section}
{\normalfont%\sffamily
  \bfseries%\color{cyan}
}
{\thesection}{1ex}{}
\titleformat{\subsection}
{\normalfont%\sffamily
  \bfseries%\color{cyan}
}
    {\thesubsection}{1ex}{}
\newcommand{\Qed}{\hspace*{10mm}\hfill\rule[-2pt]{6pt}{9pt}}
\newcommand{\W}{{\bb W}_{2}^{1}\hspace{-.62cm}{~}^{{~}^{{~}^{o}}}\,\,\,}
\newcommand{\NN}{\ensuremath{\mathbb{N}}}
\newcommand{\LL}{\widetilde{\boldsymbol{A}}} 
\newcommand{\LAL}{\widetilde{\boldsymbol{\mathfrak{A}}}}
\newcommand{\LCL}{\widetilde{\boldsymbol{\mathfrak{C}}}}
\newcommand{\LIL}{\widetilde{\boldsymbol{\mathfrak{I}}}}
\newcommand{\LSL}{\widetilde{\boldsymbol{\mathfrak{S}}}}
\newcommand{\PP}{\Upsilon}
\newcommand{\di}{\mathrm{div}\,}

\newcommand{\Ra}{\mathrm{Ra}}
\newcommand{\pti}{\partial_t}
\renewcommand{\Pr}{\mathrm{Pr}}

\newcommand{\bu}{\underline{{u}}}
\newcommand{\bv}{\underline{{v}}}
\newcommand{\bw}{\underline{{w}}}
\newcommand{\bz}{\underline{{z}}}
\newcommand{\cA}{{\mathcal{A}}}

\newcommand{\bx}{\underline{x}}
\newcommand{\bephi}{{\underline{\mathfrak{e}}}_{{\;\!}\varphi}}
\newcommand{\ber}{{\underline{\mathfrak{e}}}_{{\;\!}r}}
\newcommand{\bedrei}{{\underline{\mathfrak{e}}}_{{\;\!}{3}}}
\newcommand{\beeis}{{\underline{\mathfrak{e}}}_{{\;\!}{1}}}

\begin{document}
%%%%%%%%%%%%%%%
\hspace*{1cm}\\[2cm]
{\large \textbf{Natural convection in the horizontal annulus: critical Rayleigh
    number for the steady problem} }
\\[.5cm]
Arianna Passerini\\
\hspace*{2mm}{\small \em {Department of Mathematics and Scientific Computing, Universit{\`a} di Ferrara,
Via Saragat,~1, 44100
Ferrara, Italy}}
\\[2mm]
Bernd Rummler\\
\hspace*{2mm}{\small \em 
  Institute for Analysis and Numerics, Otto-von-Guericke-Universit\"at
  Magdeburg,  PF
  4120, 39016 Magdeburg, Germany}
\\[2mm]
Michael \rose\\
\hspace*{2mm}{\small \em Institute for Applied Mathematics,
  Universit{\"a}t Freiburg, 
  Eckerstr.~1, 79104 Freiburg, Germany}
\\[2mm]
Gudrun Th\"ater\footnote{Corresponding author\quad
  E-mail:~\textsf{gudrun.thaeter@kit.edu}}\\
\hspace*{2mm}{\small \em  Institute for Applied \& Numerical Mathematics,
  KIT, 
  76128 Karlsruhe, Germany}
\\[5mm]
%%%%%%%%%%%
\hspace*{.1cm}\hfill\parbox{16cm}{
{\small\textbf{Abstract:}
For the 2D Oberbeck-Boussinesq system in an annulus we are looking
for the critical Rayleigh number for which the (nonzero) basic flow loses 
stability. For this we consider the corresponding Euler-Lagrange
equations and construct a precise functional analytical frame for the Laplace-
and the Stokes problem as well as the Bilaplacian operator in this
domain. With this frame and the right set of basis functions it is then possible to
construct and apply a numerical scheme providing the critical Rayleigh
number. These results will be published in ZAMM in 2025.
}}
%%%%%%%%%%%
\section{Introduction}
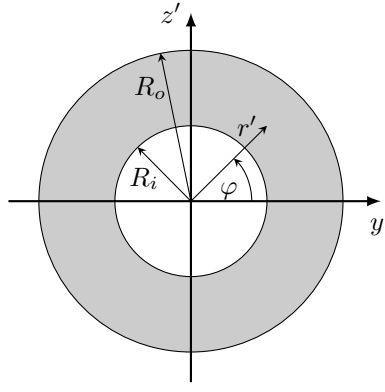
\begin{figure}[t]
\centering
  \begin{tikzpicture}
    \filldraw[fill=black!20] (0,0) circle (2);
    \filldraw[fill=white] (0,0) circle (1);
    \draw[->, >=latex, thick] (-2.4,0)--(2.5,0);
    \draw[->, >=latex, thick] (0,-2.4)--(0,2.5);
    \draw (2.5,0) node[anchor=north] {$y'$};
    \draw (0,2.5) node[anchor=east] {$z'$};
    \draw[->, >=stealth] (0,0)--(1,1);
    \draw(1,1) node[anchor=east] {$r'$};
    \draw[->, >=stealth] (0.8, 0) arc (0:45:0.8);
    \draw(0.5, -0.1) node[anchor=south] {$\varphi$};
    \draw[->, >=stealth] (0,0)--(-0.4,1.96);
    \draw(-0.2,1.5) node[anchor=east]{$R_o$};
    \draw[->, >=stealth] (0,0)--(-0.7,0.714);
    \draw(-0.3,0.3) node[anchor=east]{$R_i$};
  \end{tikzpicture}
\caption{Twodimensional flow region --
temperature boundary values $T_i>T_o$.}
\label{fig:cylindric_flow_region}
\end{figure}

%%%%%%%%%%%%%%%%
Consider the flow of an incompressible Newtonian fluid between two
horizontal coaxial cylinders with radii $0<R_i< R_o$ (cf.~Fig.~\ref{fig:cylindric_flow_region}).
The flow is driven by a
gravitational field perpendicular to 
the cylindrical axis and the 
temperature difference between $T_i$ on the inner and
$T_o$ on the outer jacket (where 
$T_o<T_i$).

This setting is a model to study very diverse phenomena, such as thermal energy storage
systems, aircraft cabin insulation, cooling of electronic components,
electrical power cable, and thin films.%Hier ergänzen?

In this configuration the flow is mainly characterized by two
nondimensional parameters:
the ``thinness'' of the gap between the cylinders, which we measure as the
{\em inverse relative gap width}
\begin{equation}\label{calA}
\mathcal{A}:=\frac{2R_i}{R_o-R_i}
\quad\mbox{and}\quad \Ra:=\frac{\alpha g}{\nu k}(T_i-T_o)(R_o-R_i)^3\,,
\end{equation}
the {\em Rayleigh number}, which classifies the heat transfer regime in the
flow. More precisely, it measures the ratio between
conduction and convection since in the definition we have
 $\alpha$ as 
the volumetric expansion coefficient, $g$ the gravity
acceleration, $\nu$ the kinematic viscosity, and $k$ stands for
the thermal diffusivity.

The study of convective flows between horizontal coaxial cylinders inside
a gravitational field 
shows (see references in \cite{FPRT}) that only for an intermediate range of the parameter 
$\cal A$ it is interesting to look
at three-dimensional behaviour, and for ${\cal A}<2,8$ (wide gap) as
well as ${\cal A}> 8,5$ (small gap) everything interesting happens in the
cross-section.
For that in this paper we study the problem in the two-dimensional annulus (i.e. the
cross-section of the above mentioned geometry) knowing that our study
will also give information for the 
three-dimensional flow problem in the small and wide gap cases.

In contrast to the classical Bénard problem in our geometry even for
very small Rayleigh numbers there is never a zero velocity
solution. Instead there is a so-called basic flow
which solves our equations for any $\Ra$ but for which we do not know an analytical
expression. For that the simplest methods which work for the Bénard problem
cannot be adapted for the annulus. But we can still work with the characterisation of
the critical Rayleigh number $\Ra_c$ as the inverse of the supremum of the
ratio between functionals 
for the kinetic
energy and the heat transfer in  \eqref{3.5} below  (see, e.g.,
\cite{Straughan}). In order to calculate our $\Ra_c$ with 
this supremum we first derive  the corresponding Euler-Lagrange equations and
later on transform them to an eigenvalue problem for a compact
self-adjoint operator. For this we have to
make sure that the abstract procedure translates into a precisely known
frame for our operators. Here our paper becomes a bit technical since
we need a few convenient
properties of all involved operators (which partially are almost
obvious and partially need a few lines of proof and functional
analytical results). In particular we prove that the
supremum exists and is indeed a maximum in Lemma
\ref{maxex}. Moreover, it is finite and thus, the critical Rayleigh
number as defined in \eqref{3.5} is finite as well.

Unfortunately, it is not
straightforward to use standard eigenvalue software in order to
numerically solve the found eigenvalue problem for the critical
Rayleigh number since it turns out not to be stable. Instead we use
special sets of bases of 
eigenfunctions of the Laplace and the Stokes problem in the annulus
(which contain Bessel functions) to formulate and apply 
a numerical scheme which is stable and approximates the wanted
$\Ra_c$ without having to refer to determinants which are too big to
handle comfortably. Here we have to translate the representation of operators in the
existence proof to a (slightly) different set of 
operators which make the numerical scheme stable.

Finally, in two-dimensional flow  people often prefer to work with the
streamfunction formulation. For
that we include how to formulate the functional
analytical frame with the streamfunction and the Bilaplace operator in
the annulus.

Let us introduce the necessary {\em  notation}.\\[.2cm]
{\bf General notation A.} Let ${\bb R}^{2}$ be endowed with the
usual Euclidian norm $\| . \|$ 
and  elements of ${\bb R}^{2}$ be denoted by underlined small letters.
The unit circle is  $\omega\,:=\,\{{\underline{x}}\in {\bb
  R}^{2}:\,\|{\underline{x}}\|=1\}$ 
and closed circles around the origin with radius $r$ are
$\omega_{r}\,:=\,\{{\underline{x}}\in {\bb
  R}^{2}:\,\|{\underline{x}}\|=r\} 
=r\omega$
for all $r\,\in\,(0,\infty)$.
\\[3mm]
{\bf Annulus domains.} It is useful to study our domain without
dimensions. As usual 
we pick the annulus with fixed gap width $1$ using our non-dimensional
parameter 
${\cal A}$ as follows:
For any ${\cal A}\,\in\,(0,\infty)$, we denote
by
\[
{\Omega}_{{\cal A}}\,:=\{{\underline{x}} \in {\bb R}^{2}:\,{\cal A}/{2}
<\|{\underline{x}}\|<1+{\cal A}/{2}\}\,.
\]
Its boundary  
$\partial{\Omega}_{\cal A}$ consists of two parts, namely the inner
and outer boundary
$\omega_{{\cal A}/{2}}$ and  $\omega_{1+{\cal A}/{2}}\,,$ respectively.
\\[3mm]
{\bf General notation B.} Let ${\Omega}$ be any of the domains
introduced above.  In what follows we will use by way of an
abbreviation $(.)$ for $({\Omega})$ everywhere. 
We consider the usual Lebesgue and Sobolev spaces ${\bb L}_{2}(.)$ and
${\bb W}_{2}^{k}(.)$ 
of scalar functions and the Lebesgue and Sobolev spaces 
${\underline{\bb L}}_{{\;\!}{2}}(.)=({\bb L}_{2}(.))^{2}$
and ${\underline{\bb W}}_{{\;\!}{2}}^{k}(.)=({\bb
  W}_{2}^{k}(.))^{2}$ of vector functions. 
The scalar product in ${\bb L}_{2}(.)$ and ${\underline{\bb L}}_{{\;\!}{2}}(.)$ is written as $(.,.)_{2}\,:=\,(.,.)$ and the norms are denoted by $\| . \|_{2}$.
The notation ${\bb W}_{2}^{k}\hspace{-.62cm}{~}^{{~}^{{~}^{o}}}\hspace{.2cm}(.)$ 
is taken for the closure of $C_{o}^{\infty}(.)$ 
in ${\bb W}_{2}^{k}(.)$.
All solenoidal vector functions belonging to 
${\underline{C}}_{{\;\!}{o}}^{\infty}(.)$ form $\underline{\cal V}(.)$. The closures
of $\underline{\cal V}(.)$ in ${\underline{\bb L}}_{{\;\!}{2}}(.)$ 
and
${\underline{\bb W}}_{{\;\!}{2}}^{1}(.)$, respectively, are denoted by
${\underline{\bb H}}(.)$ 
and ${\underline{\bb V}}(.)$, respectively. We suppose for all function written in polar
coordinates $(r,{\varphi})$ the general periodicity in the angular coordinate 
${\varphi}$. We define the dual spaces 
%explain???was soll das heissen? 
%the spaces of Distribution
   {{${\bb W}_{2}^{-2}(.):=({\bb
      W}_{2}^{2}\hspace{-.62cm}{~}^{{~}^{{~}^{o}}}\hspace{.2cm}(.))'$}},
${\bb W}_{2}^{-1}(.):=({\bb
  W}_{2}^{1}\hspace{-.62cm}{~}^{{~}^{{~}^{o}}}\hspace{.2cm}(.))'$ and
${\underline{\bb V}}'(.):=({\underline{\bb V}}(.))'$ and use the
 evolution (Gelfand) triples 
\begin{equation}\label{tripels} {{[{\bb W}_{2}^{2}\hspace{-.62cm}{~}^{{~}^{{~}^{o}}}\hspace{.2cm}(.),{\bb L}_{{\;\!}{2}}(.),
{\bb W}_{2}^{-2}(.)]}}\,\,,\quad
[{\bb W}_{2}^{1}\hspace{-.62cm}{~}^{{~}^{{~}^{o}}}\hspace{.2cm}(.),{\bb L}_{{\;\!}{2}}(.),
{\bb W}_{2}^{-1}(.)]\,\quad \text{and} \quad \,[{\underline{\bb V}}(.),{\underline{\bb H}}(.),{\underline{\bb V}}'(.)]\,.
\end{equation}
%%%%%%%%
{\bf General notation C.} Due to the shape of our domain together with Cartesian we also use
polar coordinates - whichever makes more sense. In particular, $\ber$
is the unit vector in direction $r$,
$\bedrei=\sin\varphi\,\ber+\cos\varphi\,\bephi$ the unit vector in direction of
$z$, and we collect all arising gradients in $\Pi$ (thus, it is not
the thermodynamical pressure).
\begin{equation}\label{1.5}
b :=\ln \frac{R_o}{R_i}=\ln\left(1+\frac{2}{\mathcal{A}}\right)\,,
\end{equation}
is a purely geometric parameter, unbounded as $\mathcal{A}$ tends to
zero. Moreover, we have a 
number for material properties, the Prandtl number $\Pr=\nu/k$.
Instead of the temperature $T$ we treat the {\em excess temperature}
\begin{equation}%\label{1.6}%6
\tau:=\frac{T}{T_i-T_o}-T^*=\frac{T}{T_i-T_o}-\frac{T_i}{T_i-T_o}+\frac{1}{b}\left(\ln r -\ln
\frac{R_i}{R_o-R_i}\right)=\frac{T}{T_i-T_o}-\frac{T_i}{T_i-T_o}+\frac{1}{b}\ln\frac{2r}{\mathcal{A}}
\,.
\end{equation}
Here, the scalar field $T^*$ is the conductive
solution - namely,
$T^*$ solves $\, \Delta T^*=0\, $ with boundary conditions $\,
T^*(\frac{\mathcal{A}}{2},\varphi)=\frac{T_i}{T_i-T_o}\, $ and $\,
T^*(\frac{\mathcal{A}}{2}+1,\varphi)=\frac{T_o}{T_i-T_o}\, $.\\[1mm]
{\bf Equations.} The full non-dimensional Oberbeck-Boussinesq system in polar
coordinates on $(0,\infty)\times\Omega_{\mathcal{A}}$ is (see, e.g., \cite{FPRT})
\begin{equation}\label{oursystem}
\begin{aligned}
\frac{1}{\Pr}(\pti \bv +({\underline{\nabla}}^T\cdot\bv)^T\cdot\bv)-\Delta\bv+{\underline{\nabla}}\Pi
&=
\dfrac{\Ra}{{b}}\sin\varphi\,\ber+\Ra\,\tau\bedrei\,,
\\ 
\di\bv&=0\,,
\\
\pti \tau+\bv^T\cdot{\underline{\nabla}}\tau -\Delta\tau
&=\frac{1}{r\,{b}} \bv^T\cdot \ber
\end{aligned}
\end{equation}
endowed with the boundary and initial conditions
\begin{equation}\label{bc}
\bv={\underline{0}}\,,\quad
\tau=0\quad\mbox{on}\quad\partial\Omega_{\mathcal{A}}
\qquad\mbox{and}\qquad
\bv=\bv_0\,,\quad
\tau=\tau_0 \quad\mbox{for }t=0\,.
\end{equation}
%%%%%%%%%%%%%%%%%%%%%%%%%%%%%%%%%
In \cite{FPRT} it is discussed in detail that the stability analysis of
\eqref{oursystem} is closely related to the investigation of the
stability of the basic flow which solves the following linear 
homogeneous system 
\begin{align}\label{3.4}
{\underline{\nabla}}^T \cdot\bw&=0\,, \nonumber\\
\frac{1}{\Pr}{\partial_t\bw}
-\Delta\bw +{\underline{\nabla}} p&= \Ra\,\theta\,\bedrei\,,\\
{\partial_t\theta}-\Delta\theta &= \frac{
w_r}{r\, b}\nonumber\,
\end{align}
with $\bw$ and $\theta$ as unknowns and zero Dirichlet boundary
conditions. In particular, it is shown in \cite[Sec.~4]{FPRT} that
the number $\Ra_c$ (defined below) is a good approximation of the critical
Raleigh number for the asymptotic non-linear stability of steady flows
for small $\mathcal A $.  We study (\ref{3.4}) in detail.

The physical consistency of the Oberbeck-Boussinesq system is not
discussed here, but we want to point out that although the isochoric
condition is in the system, it can hold true also for
compressible fluids, such as perfect gases, as recently shown in
\cite{GP}. Compared to the classical Oberbeck-Boussinesq system, the system rigorously
derived in \cite{GP} is different only for an extra independent
nondimensional parameter.

\section{{{Description by functionals and Euler-Lagrange equations}}}
%%%%%%%%%%%%%%%%%%%%%%%%%%%%%%%%%%%%%%%%%%%%%%%%%%%%%%%%%%%%%%%
The validity of {\em  Poincar\'e inequalities} 
ensures, that  the spaces ${\bb
  W}_{2}^{1}\hspace{-.62cm}{~}^{{~}^{{~}^{o}}}\hspace{.2cm}(.)$ and  
${\underline{\bb V}}(.)$ can be equipped with equivalent
norms - which otherwise would only be semi-norms. Namely, we introduce
the Dirichlet scalar products and Dirichlet norms using Cartesian
coordinates as follows: 
\begin{align} \label{Dirichlet} 
  \begin{aligned}
    (u,v)_{D} := \,\sum_{k=1}^{2}( {\frac{\displaystyle{\partial
          u}}{\displaystyle{\partial x_{k}}}},
    {\frac{\displaystyle{\partial v}}{\displaystyle{\partial
          x_{k}}}})_2\,\, \,\forall \, u, v\,\in\,{\bb
      W}_{2}^{1}\hspace{-.62cm}{~}^{{~}^{{~}^{o}}}\hspace{.2cm}(.)\,,\quad
    (\bu,\bv)_{D} := \,\sum_{j,k=1}^{2}( {\frac{\displaystyle{\partial
          u_{j}}}{\displaystyle{\partial x_{k}}}},
    {\frac{\displaystyle{\partial v_{j}}}{\displaystyle{\partial
          x_{k}}}})_2\quad \forall \, \bu, \bv \,\in\, {\underline{\bb
        V}}(.)\,
    % \nonumber
    \\
    \|u\|_{D} :=\, \sqrt{(u,u)_{D} } \quad \forall \, u\,\in\,{\bb
      W}_{2}^{1}\hspace{-.62cm}{~}^{{~}^{{~}^{o}}}\hspace{.2cm}(.)\,,
    \quad \quad
    \|{\underline{u}}\|_{D} := \,\sqrt{ (\bu,\bu)_{D} } \quad
    \quad\quad \quad \forall \, {\underline{u}} \,\in\,{\underline{\bb
        V}}(.)\,,
  \end{aligned}
\end{align}
where the so-called Frobenius inner product is involved in the
definition.
\\[.2cm]
In this paper we provide a procedure to compute the critical Rayleigh number ${\Ra}_c$
starting with a standard method, e.g., from Straughan \cite{Straughan}.
So, for fixed $\mathcal{A}$ and given functions
$(\bw ,\theta)\in {\underline{\bb V}}\times \W$, we define 
\begin{align}
  \begin{split}
    D(\bw, \theta)&:=\| \bw\|_D^2+\|\theta\|_D^2\, ,\\
  \mathcal{F}(\bw, \theta)&:=
  (\theta,{w}_z)+\frac{1}{b}\left(\theta,\frac{w_r}{r}\right)\,
  ,\qquad \mbox{where }w_z:=\bw^T\cdot\bedrei\,,
  \end{split}\label{Ddef}
  \\\label{3.5}
  \dfrac{1}{{\Ra}_c(\cA)}&:=\sup\frac{\mathcal{F}(
    \bw,\theta)}{D(\bw, \theta)}\, .
\end{align}
The supremum in \eqref{3.5} is calculated for all non-trivial  couples of functions
$(\bw,\theta)\in {\underline{\bb V}}\times \W$.
\begin{lem}
  \label{maxex}
  The supremum in \eqref{3.5} is attained as maximum by a pair of functions
  fulfilling the Euler-Lagrange equations. This means in particular
  that a critical point
  of $F/D$ exists.
\end{lem}
\noindent{\bf Proof.} Let ${\underline{\mathfrak{w}}}:=(\bw,\theta)\,\in\,{\underline{\bb V}}\times \W:=X$. We define
\begin{equation*}
J(\underline{\mathfrak{w}}):=\frac{\mathcal{F}(\bw,\theta)}{D(\bw, \theta)}\quad\mbox{and}\quad  J_0 :=\sup J\,,
\end{equation*}
where the supremum is taken as in \eqref{3.5}.
We can directly check that $J$ is homogeneous of degree zero in real
$\lambda\neq 0$, or explicitly  
$J(\underline{\mathfrak{w}})=J(\lambda
\underline{\mathfrak{w}})$  $\forall\,\lambda\neq 0$
and thus values of $J$ are 
determined by values on the unit sphere
$\|\underline{\mathfrak{w}}\|_X = \sqrt{D(\bw, \theta)} = 1$. 
$J$ is not
infinitely large, thus $J_0<\infty$. To see this, we use the
boundedness of $\mathcal{F}(.,.)$, eq. \eqref{Ddef} for  
$\underline{\mathfrak{w}}:=(\bw,\theta)\,\in\,{\underline{\bb
    H}}\times {\bb L}_{{\;\!}{2}}:=X_1$ and the 
compact embedding of $X$ in $X_1$ 	
as well as $\|\underline{\mathfrak{w}}\|_X = 1$. Therefore the set
$\mathcal{F}(\{\underline{\mathfrak{w}}: 
\|\underline{\mathfrak{w}}\|_X = 1
\})$ is compact. Due
to compactness we observe 
$\lim \mathcal{F} (\underline{\mathfrak{w}}_n)= \mathcal{F}
(\underline{\mathfrak{w}}_0)$ and it is easy to see that $\mathcal{F}
(\underline{\mathfrak{w}}_0)$ is bounded from above by a constant. For
that there exists a sequence
$\{\underline{\mathfrak{w}}_n\}_{n=1}^\infty$  ($\mathfrak{w}_n\in
X$ for all elements of the sequence)
for which $J(\underline{\mathfrak{w}}_n)$ converges to $J_0$. Its projection to the unit sphere
converges to $J_0$ as well. Then there is a subsequence which converges
weakly in $X$ to some value $\underline{\mathfrak{w}}_0\in X$,  and
 this subsequence
strongly in $X_1 = {\underline{\bb H}}\times {\bb L}_2$ and almost everywhere
$J(\underline{\mathfrak{w}}_n)\to J(\underline{\mathfrak{w}}_0)$ and 
thus, the supremum in \eqref{3.5} exists. Next
we check due to lower semicontinuity and since $\mathcal{D}
(\underline{\mathfrak{w}}_n)$ exists, that
\begin{equation*}
  J_0=\lim  J(\underline{\mathfrak{w}}_n)=\frac {\lim \mathcal{F} (\underline{\mathfrak{w}}_n)}{\lim \mathcal{D} (\underline{\mathfrak{w}}_n)}\le \frac
  {\mathcal{F} (\underline{\mathfrak{w}}_0)}{\mathcal{D} (\underline{\mathfrak{w}}_0)}\le J_0
 \qquad\Rightarrow\quad J(\underline{\mathfrak{w}}_0)=\sup_{\underline{\mathfrak{w}}\in X}J(\underline{\mathfrak{w}})\,.   \end{equation*}
If the supremum in \eqref{3.5} is attained for some\footnote{The notation
  {{$(\widetilde{\bw},\widetilde{\theta})$}} with tilde means that the
  stationary point has nothing to do with the solutions of
  (\ref{3.4}); as we will see, it is solution of a different system of
  equations.}
$(\widetilde{\bw},\widetilde{\theta})$, then $(\widetilde{\bw},\widetilde{\theta})$  solve the
corresponding Euler-Lagrange equations, which we derive now. For an arbitrary 
real parameter $\eta$ and arbitrary {\em variations}
$(\widetilde{\bu},\widetilde{\sigma})\in {\underline{\bb V}}\times
\W$ we consider the real function
\[
  \frac{\mathcal{F}}{D}=\frac{\mathcal{F}(\widetilde{\bw}+\eta \widetilde{\bu},\widetilde{\theta}+\eta\widetilde{\sigma})}
 {D(\widetilde{\bw}+\eta
  \widetilde{\bu},\widetilde{\theta}+\eta\widetilde{\sigma})}\,.
\]
Since $(\widetilde{\bw},\widetilde{\theta})$ is a maximum the Gateaux
derivative $\frac {d}{d\eta}|_{\eta=0}$ has to vanish. Thus, 
  \begin{align}
0=\left .\frac {d}{d\eta}\left(\frac{\mathcal{F}}{D}\right)\right|_ {\eta=0}=
\left [\frac{1}{{D}^2}\left(D\frac {d}{d\eta}\mathcal{F}-\mathcal{F}\frac {d}{d\eta}
D\right)\right]_ {\eta=0}=\left .\frac{\frac {d\mathcal{F}}{d\eta}}{D}\right|_
{\eta=0} -\frac{1}{{\Ra}_c}\left .\frac{\frac {dD}{d\eta} }{D}\right|_
{\eta=0}\,.\label{eq:el}
  \end{align}
Inserting the expressions 
\begin{align*}
\left.\frac {d\mathcal{F}}{d\eta}\right|_{\eta=0}&=(\widetilde{\sigma},\widetilde{w}_z)+(\widetilde{\theta},\widetilde{u}_z)
+\frac{1}{b}\Big(\big(\widetilde{\theta},
\frac{\widetilde{u}_r}{r}\big)+\big(\widetilde{\sigma},\frac{\widetilde{w}_r}{r}\big)\Big)\,,\\
\left.\frac {dD}{d\eta}\right|_{\eta=0}&=2(\widetilde{{\bu}},\widetilde{\bw})_D+2(\widetilde{\theta},
\widetilde{\sigma})_D
\end{align*}
into \eqref{eq:el}, we conclude that for all variations
$(\widetilde{\bu},\widetilde{\sigma})\in {\underline{\bb V}}\times
\W$ it holds
\begin{align}\label{bilinear_formulation} 
(\widetilde{\sigma},\widetilde{w}_z)+(\widetilde{\theta},\widetilde{u}_z)
+\frac{1}{b}\Big(\big(\widetilde{\theta},
\frac{\widetilde{u}_r}{r}\big)+\big(\widetilde{\sigma},\frac{\widetilde{w}_r}{r}\big)\Big)
=\frac{2}{{\Ra}_c}\Big((\widetilde{{\bu}},\widetilde{\bw})_D+(\widetilde{\theta},
\widetilde{\sigma})_D\Big)\,. \Qed
\end{align}
If we assume higher regularity of
$(\widetilde{\bw},\widetilde{\theta})$ we can apply integration by
parts to find
\[
(\widetilde{\sigma},\widetilde{w}_z)+(\widetilde{\theta},\widetilde{u}_z)
+\frac{1}{b}\Big(\big(\widetilde{\theta},
\frac{\widetilde{u}_{r}}{r}\big)+\big(\widetilde{\sigma},\frac{\widetilde{w}_r}{r}\big)\Big)=
-\frac{2}{{\Ra}_c}\Big((\Delta\widetilde{\bw},\widetilde{\bu})+(\Delta\widetilde{\theta},\widetilde{\sigma})\Big)\,
.
\]
We have the freedom to choose either a) arbitrary $\widetilde{{\bu}}$ and
$\widetilde{\sigma}=0$; or b) $\widetilde{{\bu}}=0$ and
arbitrary $\widetilde{\sigma}$.
In the first case, we obtain
\[\int_{\Omega_\mathcal{A}}\big(\widetilde{\theta}\bedrei+\frac{\widetilde{\theta}}{br}\ber+\frac{2}{{\Ra}_c}
\Delta\widetilde{\bw}\big)^T\cdot\widetilde{\bu}\,d\:\!\bx =0\qquad
\Rightarrow\qquad
\widetilde{\theta}\bedrei+\frac{\widetilde{\theta}}{br}\ber+\frac{2}{{\Ra}_c}\Delta\widetilde{\bw}
+{\underline{\nabla}}\,\widetilde{p}=0\, .\]
Here we used that  $\widetilde{\bu}$ is
arbitrary but divergence free, which implies that the term in brackets belongs
to the orthogonal complement with respect to Helmholtz's decomposition. That's
why ${\underline{\nabla}}\,\widetilde{p}$ is added.
In the second case we get
\[
\int_{\Omega_\mathcal{A}}\left(\widetilde{w}_z+\frac{\widetilde{w}_r}{br}+\frac{2}{{\Ra}_c}\Delta\widetilde{\theta}\right)
\widetilde{\sigma}d\:\!\bx=0\qquad\Rightarrow\qquad
\widetilde{w}_z+\frac{\widetilde{w}_r}{br}+\frac{2}{{\Ra}_c}\Delta\widetilde{\theta}=0\,
.
\]
Moreover, if we set
\begin{align}\label{gradS}
\bedrei+\frac{\ber}{br} =\,{\underline{\nabla}}\:\! S,\qquad
\text{where}\quad S:= r\sin\varphi  +\frac{1}{b}\ln r,
\end{align}
 the strong form of the 
Euler-Lagrange equations for the maximum
$(\widetilde{\bw},\widetilde{\theta})$  of the functional $\mathcal
F/D$ read
\begin{align}\label{3.18}
\widetilde{\theta}{\underline{\nabla}}\:\! S +\lambda\Delta\widetilde{\bw}&=
-{\underline{\nabla}}\,\widetilde{p}\,,\quad {\underline{\nabla}}^T
                                                                            \cdot\widetilde{\bw}=0\,,\\
  \widetilde{\bw}^T\cdot {\underline{\nabla}}\:\! S  
\,+\lambda\Delta\widetilde{\theta}&=0 
\qquad\mbox{with}\qquad \lambda=\frac{2}{{\Ra}_c}\,.\label{3.19}
\end{align}
{{In what follows we regard the problem (\ref{3.18})-(\ref{3.19}) in the form of\\[.15cm]
{\bf Problem A.}  {\em (Euler-Lagrange equations of the functional $\mathcal
F/D$)}\\
We look for nontrivial solutions $(\lambda,{\bw},\theta)\in {\bb C}\times {\underline{\bb V}}\times \W$ of the equations
\begin{align}\label{E1}
{\theta}{\underline{\nabla}}\:\! S +\lambda\Delta {\bw}&=
-{\underline{\nabla}}\,{p}\,,\quad {\underline{\nabla}}^T
                                                                            \cdot\widetilde{\bw}=0\,,\\ 
{\bw}^T\cdot {\underline{\nabla}}\:\! S
\,+\lambda\Delta {\theta}&=0
\qquad\mbox{with}\qquad \lambda\,\in\,{\bb C}\,,\label{E2}\
\end{align}
where (\ref{E1}), (\ref{E2}) are equations in ${\underline{\bb V}}'(.)$ 
resp. ${\bb W}_{2}^{-1}(.)$.}}\\[3mm]
We will show that Problem A can be equivalently formulated as an
eigenvalue problem, which possesses a solution for all non-trivial eigenvalues $\lambda$
with the corresponding eigenvectors. We then define the critical
Rayleigh number $\Ra_c$ via the largest eigenvalue $\lambda_{\max} $ through
${\Ra}_c:=\frac{2}{\lambda_{\max}}$.
We now prepare the corresponding framework.
%%%%%%%%%%%%%%%%%%%%%%%%%%%%%%%%%%%%%%%%T3.5
\section{Standard differential operators and the symmetric operator
      $\LL$} 
\label{SecStandOps}
%%%%%%%%%%%%
In what follows we take  ${\Omega}:={\Omega}_\mathcal{A}$ as symbol for any of the annulus domains.
We will introduce the Laplacian, 
the Stokes operator, and the Bilaplacian 
in the sense of functional analysis. The domains of definition of these operators
are dense in the Hilbert spaces ${\bb L}_{{\;\!}{2}}({\Omega})$, 
${\underline{\bb H}}({\Omega})$
and ${\bb L}_{{\;\!}{2}}({\Omega})$, respectively.\\[.1cm]
{\bf Definition A.}  {\em The Laplace operator is defined as}
\begin{align*}
{\boldsymbol L^{\circ}}\,{v} := 
- \Delta{\;\!}_{\underline{x} } v
\hspace*{0.7cm}  \forall\,v\,\in\,D({\boldsymbol
  L^{\circ}})=C_{o}^{\infty}({\Omega})\,. 
\end{align*}
We denote the Friedrichs' extension of ${\boldsymbol L^{\circ}}$
by ${\boldsymbol 
  L}:={\overline{\boldsymbol L^{\circ}}}$, where ${\boldsymbol
  L}$ is defined on 
$D({\boldsymbol L})\,:=\,{\bb
  W}_{2}^{1}\hspace{-.62cm}{~}^{{~}^{{~}^{o}}}\hspace{.2cm}({\Omega}) 
\cap {\bb W}_{2}^{2}({\Omega})$.\\[.15cm]
{\bf Remark:} The range of ${\boldsymbol L}$ is
$R({\boldsymbol L})={\bb L}_{2}({\Omega})$. In this sense we write:
${\boldsymbol L}=-\Delta{\;\!}_{\underline{x} } : D({\boldsymbol
  L})\subseteq {\bb L}_{2}(\Omega)\longmapsto \,{\bb L}_{2}(\Omega)$.
\\[.1cm] 
We need some preparations to define the Stokes operator. 
The Leray-Helmholtz projector $\Upsilon$ is the well-defined
projector of ${\underline{\bb L}}_{{\;\!}2}(\Omega)$ onto its subspace 
${\underline{\bb H}}(\Omega)$
of generalised solenoidal fields with vanishing generalised traces in
the normal direction on the boundary. 
We note, that the Leray-Helmholtz projector $\Upsilon$
is also used in the sense of:
$\Upsilon :{\underline{\bb W}}_{{\;\!}2}^{1}(\Omega)\,\longmapsto
\,{\underline{\bb V}}(\Omega)$\,.\\[.1cm]
{\bf Definition B.} 
{\em The Stokes operator is defined as}
$
{\boldsymbol S^{\circ}}\,{\underline{v}}:=-
\Delta_{\underline{x} } {\underline{v}}\hspace*{0.5cm}
\forall\,{\underline{v}}\in D({\boldsymbol
  S^{\circ}})=\underline{\cal V}({\Omega})\,. 
$
We denote the Friedrichs' extension of ${\boldsymbol S^{\circ}}$ by 
${\boldsymbol S}:={\overline{\boldsymbol S^{\circ}}}$, where
${\boldsymbol S}$ is defined on
$D({\boldsymbol S})={\underline{\bb S}}_{{\;\!}}^{2}(\Omega)\,=\,{\underline{\bb V}}^{2}(\Omega)\,:=
{\underline{\bb W}}_{{\;\!}2}^{2}(\Omega)\cap{\underline{\bb V}}(\Omega)$ \,. \\[.1cm]
{\bf Remark:} The range of ${\boldsymbol S}$ is $R({\boldsymbol
  S})={\underline{\bb H}}(\Omega) $. In this context we use
${\boldsymbol S} =-\Upsilon \Delta{\;\!}_{\underline{x}
}:{\underline{\bb S}}_{{\;\!}}^{2}(.) \subseteq {\underline{\bb H}}(\Omega)
\longmapsto \,{\underline{\bb H}}(\Omega)$.
\\[.2cm]
{{Finally we define the Bilaplacian:\\[.15cm]
{\bf Definition C.} 
{\em The Bilaplacian (biharmonic operator) is defined as}
\begin{align*}
{\boldsymbol B^{\circ}}\,{v} := 
 \Delta{\;\!}^2_{\underline{x} } v
\hspace*{0.7cm}  \forall\,v\,\in\,D({\boldsymbol
  B^{\circ}})=C_{o}^{\infty}({\Omega})\,. 
\end{align*}
We denote the Friedrichs' extension of ${\boldsymbol B^{\circ}}$
by ${\boldsymbol 
  B}:={\overline{\boldsymbol B^{\circ}}}$, where 
$D({\boldsymbol B})\,:=\,{\bb
  W}_{2}^{2}\hspace{-.62cm}{~}^{{~}^{{~}^{o}}}\hspace{.2cm}({\Omega}) 
\cap {\bb W}_{2}^{4}({\Omega})$ is its domain.\\[.15cm]
{\bf Remark:} The range of ${\boldsymbol B}$ is
$R({\boldsymbol B})={\bb L}_{2}({\Omega})$. We write here:
${\boldsymbol B}=\,\Delta{\;\!}^2_{\underline{x} } : D({\boldsymbol
  B})\subseteq {\bb L}_{2}({\Omega})\longmapsto \,{\bb L}_{2}(\Omega)$
as well.
\\[.2cm] 
%%%%%%%%%%%
We recall the essential properties of the operators ${\boldsymbol L}$,
${\boldsymbol S}$, and ${\boldsymbol B}$ using the Stokes operator ${\boldsymbol S}$
as example:}}\\[-.35cm] 
\begin{theor}\label{thmstok}
The Stokes operator ${\boldsymbol S}$ is positive and self-adjoint.
Its inverse ${\boldsymbol S}^{-1}$ is injective, self-adjoint and compact.
\end{theor}
\noindent The proof of Theorem \ref{thmstok} is a simple modification of 
Theorems 4.3 and 4.4 in \cite{CF}. The essential tools are
the Rellich theorem and Lax-Milgram lemma.
The well-known theorem of Hilbert and
regularity results like \cite[Prop.~I.2.2]{Temam}
lead to  
\begin{lem}
\label{STOeiFU}
The Stokes operator is an operator with a pure point spectrum.
All eigenvalues $\lambda_{j}^{*}\,=\,(\kappa_{j})^2$
of  ${\boldsymbol S}$ are real and of finite multiplicity.
The associated eigenfunctions 
$\{{\underline{v}}_{j}({\underline{x}})\}_{j=1}^{\infty}$
of the Stokes operator ${\boldsymbol S}$  (counted in multiplicity)
are an orthogonal basis of
${\underline{\bb H}}(.)$ and ${\underline{\bb V}}(.)$. We obtain that
\begin{align*}
{{(a)}}&\quad {\boldsymbol 
S}{\underline{v}}_{j}=\lambda_{j}^{*}{\underline{v}}_{j}\quad\mbox{for
}\quad{\underline{v}}_{j}\in D({\boldsymbol S})
\quad\forall\,j=1,2,\dots\,;\\
{(b)}& \quad
0\,<\lambda_{1}^{*}\leq\,\lambda_{2}^{*}\,\leq\cdots\leq\,\lambda_{j}^{*}
\,\leq\cdots\quad\mbox{and}\quad 
\lim_{j\rightarrow\infty}\lambda_{j}^{*}=\infty\,;
\\[-2mm]
{(c)}&\quad
\|{\underline{v}}_{j}\|_{{\underline{\bb 
H}}}\,=1\quad\forall \, j=1,2,\dots\,.
\end{align*}
\end{lem}
\begin{notation} \label{eigenpairs}
We write the eigenpair-systems of the Laplacian ${\boldsymbol L}$, the
Stokes operator ${\boldsymbol S}$, and of the Bilaplacian ${\boldsymbol B}$
as:
\begin{equation} \label{eigenpairs2}
\,\,\{(\omega_{j})^2,\chi_{j}\}_{j=1}^{\infty}\, 
\,\qquad \text{\,,\,}\,\,\qquad
\{(\kappa_{j})^2,{\underline{v}}_{j}\}_{j=1}^{\infty}\,=\,
\{\lambda_{j}^{*},{\underline{v}}_{j}\}_{j=1}^{\infty}\,,\qquad \text{and}\,\,\qquad
\{\mu_{j}^2,\psi_{j}\}_{j=1}^{\infty}\,.
\end{equation}
We choose the systems of eigenfunctions of ${\boldsymbol L}$: $\{\chi_{j}\}_{j=1}^{\infty}$,  
of ${\boldsymbol S}$:
$\{\underline{v}_{j}\}_{j=1}^{\infty}$ and of
${\boldsymbol B}$: $\{\chi_{j}\}_{j=1}^{\infty}$ ordered by increasing eigenvalues 
and counted in multiplicity as an orthonormal basis in each case, such that:  
\begin{align*}
{\bb L}_{2}\,=\,{\overline{\text{span}
\{\chi_{j}\}_{j=1}^{\infty}}}^{\,{\bb L}_{2}}\,\qquad \text{\,,\,}\,\,\qquad\,{\underline{{\bb H}}}=\,{\overline{\text{span}
\{\bv_{j}\}_{j=1}^{\infty}}}^{\,{\underline{\bb H}}}
\,,\qquad \text{and}\,\,\qquad {\bb L}_{2}\,=\,
{\overline{\text{span}
\{\psi_{j}\}_{j=1}^{\infty}}}^{\,{\bb L}_{2}}\,.
\end{align*}
\end{notation}
\noindent {\bf Remark:} 
We refer to \cite{Rummler} for formulas defining the
complete sets of Laplace and  
Stokes eigenfunctions on  circular annuli 
\[
{\Omega}^{*}_{\sigma}\,:=\{{\underline{x}} \in {\bb R}^{2}:\,0<{\sigma}
<\|{\underline{x}}\|<1\}\,.
\] 
The (obvious) transformation rules from $\Omega^{*}_{\sigma}$ to
${\Omega}_\mathcal{A}$ are given in 
\cite{RRTZAMM}.
\begin{notation} \label{SymmOp} Let 
\begin{gather*}
\LL:{\underline{\bb V}}^2(\Omega)\times ({\bb
  W}_{2}^{1}\hspace{-.62cm}{~}^{{~}^{{~}^{o}}}\hspace{.2cm}({\Omega}) 
\cap {\bb W}_{2}^{2}({\Omega}) )\,=\,D({\boldsymbol S})\times D({\boldsymbol L})\subseteq
{\underline{\bb H}}({\Omega}) \times{\bb L}_{2}({\Omega}) 
 \longrightarrow {\underline{\bb H}}({\Omega}) \times{\bb L}_{2}({\Omega}) 
%\]
%\[
 \\
 \LL =(A_1,A_2)\quad\text{with}\quad \left\{\begin{array}{l}
A_1:{\underline{\bb V}}^2(.)\times {\bb
  W}_{2}^{1}\hspace{-.62cm}{~}^{{~}^{{~}^{o}}}\hspace{.2cm}({\Omega}) 
\cap {\bb W}_{2}^{2}({\Omega}) \subseteq
{\underline{\bb H}}({\Omega}) \times{\bb L}_{2}({\Omega}) \longrightarrow {\underline{\bb H}}({\Omega}) \\
A_2:
{\underline{\bb V}}^2(.)\times {\bb
  W}_{2}^{1}\hspace{-.62cm}{~}^{{~}^{{~}^{o}}}\hspace{.2cm}({\Omega}) 
\cap {\bb W}_{2}^{2}({\Omega}) \subseteq
{\underline{\bb H}}({\Omega}) \times{\bb L}_{2}({\Omega})  \longrightarrow {\bb L}_{2}({\Omega}) 
                                            \end{array} \right.
                                        \\
A^{\underline{\bb V}}_i:= \left.A_i\right|_{D({\boldsymbol S}) \times\{{0}\}} ;\quad A^{\bb W}_i:= \left.A_i\right|_{\{\underline{0}\}\times D({\boldsymbol L})
}\,\,\,,\quad i=1,2\,,
\end{gather*}
be the operators defined by density
$\forall\, {\underline{\Psi}}\,\in\,D({\boldsymbol S})$
and $\forall\, {\varphi}\,\in\,D({\boldsymbol L})$ via
\begin{align*}
A_i&=A_i^{\underline{\bb V}}+A^{{\bb W}}_i, \quad i=1,2\,:&\,
\\
({\underline{\Psi}}, A_1^{\underline{\bb V}}\,\bw+A_1^{{\bb W}}\theta )&= \,-\,{\lambda}
\left({\underline{\Psi}},\bw\right)_{D}\,+\,\left({\underline{\Psi}},\theta
{\underline{\nabla}}\;\!S \right) &\,{~}\, 
\\
(\varphi, A_2^{\underline{\bb V}}\,\bw+A_2^{{\bb W}}\theta )&= \,-\,{\lambda}\left(\varphi,\theta\right)_{D}
\,+\,\left(\varphi,({\underline{\nabla}}\;\!S)^T \bw\right)&\,{~}\, 
\end{align*}
 and
 \begin{equation}\label{SymmDefi}
 \LL ((\bw,\theta))=\left(\begin{array}{cc}
 A^{\underline{\bb V}}_1&A_1^{{\bb W}}\\
 &\\
 A_2^{\underline{\bb V}}&A_2^{{\bb W}}
 \end{array}\right)\left(\begin{array}{c}
 \bw\\  \\  \theta
 \end{array}
 \right) \,=\,\left(\begin{array}{cc}
- {\lambda}\, {\boldsymbol S}& {\underline{\nabla}}\;\!S\,\\
 &\\({\underline{\nabla}}\;\!S)^T \,
&- {\lambda}\, {\boldsymbol L}
 \end{array}\right)\left(\begin{array}{c}
 \bw\\  \\  \theta
 \end{array}
 \right) .
 \end{equation}
\end{notation} 
\noindent Using $\PP$ as the projection onto the divergence free
vector space we find 
\begin{align*} \big( A_1=A_1^{\underline{\bb V}}+A_1^{{\bb W}}\big)(\bw,\theta)
=\,- {\lambda}\, {\boldsymbol S}\bw +\PP\mathcal{L}\theta \quad \text{and} \quad 
\big( A_2=A_2^{\underline{\bb V}}+A_2^{{\bb W}}\big)(\bw,\theta)=\, ({\underline{\nabla}}\;\!S)^T \bw
- {\lambda}\, {\boldsymbol L}\theta \,
\end{align*} 
(cf. Theorem \ref{thmstok} and \eqref{gradS}), where
\[\mathcal{L}(\theta )= \,\theta{\underline{\nabla}} S.\]
\begin{lem} 
$\LL$ is symmetric {{in the Hilbert space ${\underline{\bb H}}(.) \times  {\bb L}_{2}(.)$}}.
\end{lem}
\noindent {\bf Proof.} We use the definition and our equations to calculate
\begin{align*}
({\underline{\Psi}}, A_1^{\underline{\bb V}}\bw)&=\,-\,\lambda\cdot({\underline{\Psi}},
\bw)_{D}=\,
(A^{\underline{\bb V}}_1{\underline{\Psi}},\bw)\quad\,,
\\
(\varphi,A^{{\bb W}}_2\theta)&=\,-\lambda\cdot(\varphi
,\theta )_{D} \quad
=(A^{{\bb W}}_2\varphi ,\theta)\quad\,,\\
({\underline{\Psi}},A^{{\bb W}}_1\theta)&
 =  ({\underline{\Psi}},\theta{\underline{\nabla}} S)=
 (\theta,({\underline{\nabla}}\;\!S)^T 
 \cdot{\underline{\Psi}})=%\quad\,,\\
  (\theta,A^{\underline{\bb V}}_2\,{\underline{\Psi}})
\end{align*}
or
 $\left(A^{\underline{\bb V}}_2\right)=\left(A_1^{{\bb W}}\right)^T$.\Qed \\[1mm]
 Notice that the eigenvalue, which is equivalent to Problem A, is
 related to the symmetric operator given by the blocks outside the
 diagonal, while the diagonal blocks multiplied with $(-1)$ give a
 symmetric metric tensor (by definition positive definite for any
 $\lambda\,>\,0$).\\ 
 We note, that the energetic extension ${\boldsymbol S}_{en}:{\underline{\bb V}} 
 \longrightarrow {\underline{\bb V}}'$ of the Stokes operator ${\boldsymbol S}$ and
the energetic extension ${\boldsymbol L}_{en}:
{\bb W}_{2}^{1}\hspace{-.62cm}{~}^{{~}^{{~}^{o}}} \longrightarrow {\bb
  W}_{2}^{-1}$ of ${\boldsymbol
  L}$ are symmetric operators in the Hilbert spaces ${\underline{\bb V}}'$ and ${\bb
  W}_{2}^{-1}$ as well. One can use the structure of \eqref{SymmDefi} to make sure that the extension of $\LL$ on ${\underline{\bb V}}'\,\times\,{\bb
  W}_{2}^{-1} $ is symmetric.
%%%%%%%%%%%%%%%%%%%%%%%%%%%%%%%%%%%%%%%%%
\section{Investigation of Problem A}
\label{Secvelocity}
Our main theoretical result  
for Problem A is the following:
\begin{theor}\label{T3.5}
For all $\mathcal{A}$ problem A (i.e. eqs.~\eqref{E1}-\eqref{E2}) is 
equivalent to the eigenvalue 
problem of the compact self-adjoint operator\, ${\tilde{{\cal C}}}\,\in\,
 {\cal L}\hspace*{-.19cm}{{\cal L}}(l\hspace*{-.075cm}{l}_{2}\times l\hspace*{-.075cm}{l}_{2},l\hspace*{-.075cm}{l}_{2}\times l\hspace*{-.075cm}{l}_{2})$ (see below). 
All eigenvalues $\lambda$
are real
and $\lambda\,=\,0$ is not an eigenvalue of the operator
$\widetilde{\cal C}$, but is 
%as 
an accumulation point of eigenvalues,
which is the only  element in the continuous spectrum of {~}\,${\tilde{{\cal C}}}$.
The critical Rayleigh number ${\Ra}_c$ is fixed through
\[\qquad %\lambda\,=
  \frac{2}{{\Ra}_c}\,:=\, |\lambda_{max}|\,=\,\|{\tilde{{\cal C}}}\|_{{\cal L  }\hspace*{-.15cm}{{\cal L}}}\,\,.\]
The set $\{\lambda\}$ of the ${\tilde{{\cal C}}}$-eigenvalues
consists of a countably infinite number of eigenvalues with finite multiplicity.\\ 
The eigen-triples
$(\lambda,{\tilde{\mathbf{a}}},
{\tilde{\mathbf{b}}})\,\in {\bb C}\times l\hspace*{-.075cm}{l}_{2}\times l\hspace*{-.075cm}{l}_{2}$
correspond to  a solution 
$(\lambda,{\tilde \bw},{\tilde \theta})\in  
{\bb C}\times{\underline{\bb V}}\times {\W}
$
of Problem A.
\end{theor}
\noindent We will prove Theorem \ref{T3.5} by functional analytical tools. 
Problem A, i.e., the 
Euler-Lagrange equations of the functional $\mathcal F/D$, is reduced to an ``algebraic one"\, by the use
of the complete system of the eigenfunctions of the Laplacian
${\boldsymbol L}$ and the Stokes operator ${\boldsymbol S}$.  
In what follows
we are going to make some arrangements for the proof of Theorem
\ref{T3.5}. 
\begin{lem}
\label{Completeness}
The system  (cf. Notation \ref{eigenpairs})
\begin{align*}
\left\{({\underline{v}}_j,\chi_k)\right\}_{j,k \in \NN}\quad 
\end{align*}
is a complete basis for the space %herein defined, namely 
 ${\underline{\bb V}}\,\times \,{{\bb
  W}_{2}^{1}\hspace{-.62cm}{~}^{{~}^{{~}^{o}}}}$\,\, .
\end{lem}
\noindent {\bf Proof.} We apply the properties of positive
self-adjoint operators with a pure point 
spectrum (cf. Th. \ref{thmstok}, Lemma \ref{STOeiFU}).
\Qed \\[.1cm]
%%%
In what follows we are going to use the Fourier representations
for any velocity $\bw\,\in \,{\underline{\bb V}}$ \,
and any temperature ${\theta}\,\in \,\W$
\,
as expansions in the systems of eigenfunctions $\{\bv_{j}\}_{j=1}^{\infty}$ 
and $\{\chi_{j}\}_{j=1}^{\infty}$ in the sense of ${\underline{\bb 
H}}$ and ${\bb L}_{2}$, respectively (cf. Lem. \ref{STOeiFU}, \eqref{eigenpairs2}):\begin{align} \label{F-expensions1}
\bw =\,\sum_{j=1}^{\infty} c_{j}\bv_{j}\,\qquad\qquad
                                                                                     {\text{and}}\qquad\qquad {\theta}=\,\sum_{j=1}^{\infty}d_{j}\chi_{j}\,. 
\end{align}
The sequences $\{c_{j}\}_{j=1}^{\infty}$ and $\{d_{j}\}_{j=1}^{\infty}$  - as
sequences of real or complex numbers - are elements of the {\em Hilbert space of sequences}
$l\hspace*{-.075cm}{l}_{2}$ because of the properties of the orthonormal 
systems $\{\bv_{j}\}_{j=1}^{\infty}$ and $\{\chi_{j}\}_{j=1}^{\infty}$ (in ${\underline{\bb 
    H}}$ and ${\bb L}_{2}$).
It is defined by 
$l\hspace*{-.075cm}{l}_{2}\,:=\,
(l_{2},\| . \|_{2,(l)})$, with the vector space
$ l_{2}\,:=\,
\{ \,\{ a_j \}_{j=1}^{\infty}\,;\,a_j \in\, {\bb R} (\text{ resp. } 
\in\, {{\bb C}} )\,\forall \, j \in\, {\bb N}: 
\,\sum_{j=1}^{\infty}|a_j|^{2}\,<\,\infty\}$ and the norm $\|
. \|_{2,(l)} $, which is
fixed through the scalar product 
$(\{ a_j \}_{j=1}^{\infty},\{ b_j \}_{j=1}^{\infty})_{2,(l)}\,:=\,
\sum_{j=1}^{\infty}\, a_j {\DarkGr{\overline{\Black b}}}_j $ \,,\, $\forall\,
\{ a_j \}_{j=1}^{\infty} \,,\{ b_j \}_{j=1}^{\infty} \,\in\, 
l_{2}$\,:\\
\hspace*{4cm}
$\| \{ a_j \}_{j=1}^{\infty} \|_{2,(l)}\,:=\,
\sqrt{\sum_{j=1}^{\infty}\, a_j {\DarkGr{\overline{\Black a}}}_j }$\quad $\forall\,
\{ a_j \}_{j=1}^{\infty} \,\in\, 
l_{2}$\,.\\[.15cm] 
We choose subspaces of $l\hspace*{-.075cm}{l}_{2}$ regarding 
any $\bw\,\in \,{\underline{\bb V}}$ or any ${\theta}\,\in \,\W$  \,
as weighted sequences,
where the eigenpairs (cf. \eqref{eigenpairs2}) will be especially highlighted.
\begin{notation} \label{SequSpacesVW} 
{\em (Special sequence spaces of coefficients)}\\ 
Let us regard the sequences of eigenvalues from the  eigenpairs-system for the Laplacian ${\boldsymbol L}$ and 
the Stokes operator ${\boldsymbol S}$ (cf. Notation \ref{eigenpairs})
ordered by increasing eigenvalues taking into account their 
multiplicities (and the correspondent eigenfunctions):
\begin{equation} \label{eigenval_sequ}
\{(\omega_{j})^2\,=\,|\omega_{j}|^2\}_{j=1}^{\infty}\,\quad ,\quad
\{\lambda_{j}^{*}\,=\,|\kappa_{j}|^2\}_{j=1}^{\infty}
\end{equation}
and the correspondent sequences of positive roots:
\begin{align*}
\{\omega_{j}\}_{j=1}^{\infty}\,\quad\quad,\quad\quad
\{\sqrt{\lambda_{j}^{*}}\,=\,\kappa_{j}\}_{j=1}^{\infty}\,\,.
\end{align*}
We declare the spaces $l\hspace*{-.075cm}{l}^{{\bb W}}_{2 }$\,,\,
$l\hspace*{-.075cm}{l}^{\underline{\bb V}}_{2 }$ with 
respect to spectral operators in the following way: 
$l\hspace*{-.075cm}{l}^{{\bb W}}_{2 }\,:=\,
(l^{{\bb W}}_{2 },\| . \|_{2,(\bb W)})$ is defined by 
\begin{eqnarray}\label{l2W_sequ}
l^{{\bb W}}_{2 }\,:=\,
\{ \,\{ a_j \}_{j=1}^{\infty}\,\in\,l_{2 }:
\,\sum_{j=1}^{\infty}\omega_{j}^2|a_j|^{2}\,<\,\infty\} \hspace*{6cm} \nonumber \\
\,\text{with the norm and the scalar product}\, \,\forall\,
\{ a_j \}_{j=1}^{\infty} \,,\{ b_j \}_{j=1}^{\infty} \,\in\, 
l^{{\bb W}}_{2 } \,:\hspace*{2cm}\\
\| \{ a_j \}_{j=1}^{\infty} \|_{2,(\bb W)}\,:=\,\sqrt{(\{ a_j
  \}_{j=1}^{\infty},\{ a_j \}_{j=1}^{\infty})_{2,(\bb W)}}\,\quad,
  \quad\, 
(\{ a_j \}_{j=1}^{\infty},\{ b_j \}_{j=1}^{\infty})_{2,(\bb W)}\,:=\,
\sum_{j=1}^{\infty}\, \omega_{j}^2 a_j {\DarkGr{\overline{\Black b}}}_j   \,.
\nonumber
\end{eqnarray}
$l\hspace*{-.075cm}{l}^{{\bb V}}_{2 }\,:=\,
(l^{{\bb V}}_{2 },\| . \|_{2,(\bb V)})$ is 
\begin{eqnarray}\label{l2V_sequ}
l^{{\bb V}}_{2 }\,:=\,
\{ \,\{ a_j \}_{j=1}^{\infty}\,\in\,l_{2 }:
\,\sum_{j=1}^{\infty}\kappa_{j}^2|a_j|^{2}\,<\,\infty\} \hspace*{6cm} \nonumber \\
\,\text{with the norm and the scalar product}\, \,\forall\,
\{ a_j \}_{j=1}^{\infty} \,,\{ b_j \}_{j=1}^{\infty} \,\in\, 
l^{{\bb V}}_{2 } \,:\hspace*{2cm}\\
\| \{ a_j \}_{j=1}^{\infty} \|_{2,(\bb V)}\,:=\,\sqrt{(\{ a_j
  \}_{j=1}^{\infty},\{ a_j \}_{j=1}^{\infty})_{2,(\bb V)}}\,\quad,
  \quad\, 
(\{ a_j \}_{j=1}^{\infty},\{ b_j \}_{j=1}^{\infty})_{2,(\bb V)}\,:=\,
\sum_{j=1}^{\infty}\, \kappa_{j}^2 a_j {\DarkGr{\overline{\Black b}}}_j   \,.
\nonumber
\end{eqnarray}
\end{notation} 
\begin{lem}
\label{Hilbert-sequences}
Let $\{ c_j \}_{j=1}^{\infty} $ be the Fourier coefficients of
$\bw\,\in\,{\underline{\bb V}}$ with respect to the  
${\underline{\bb H}}$ orthonormal system $\{\bv_{j}\}_{j=1}^{\infty}$
and $\{ d_j \}_{j=1}^{\infty} $ be the ones for 
${\theta}\,\in\, \W$ with respect to the ${\bb L}_{2}$ orthonormal system 
$\{\chi_{j}\}_{j=1}^{\infty}$.
Then the following statements are true:
\begin{align*}
{{(i)}}& \quad \{ c_j \}_{j=1}^{\infty} \,\in\, 
l\hspace*{-.075cm}{l}^{{\bb V}}_{2 }\,\, \Leftrightarrow \,\bw\,\in\,{\underline{\bb V}}\,\,,
&&\|{\underline{w}}\|_{D}\,=\,
\| \{ c_j \}_{j=1}^{\infty} \|_{2,(\bb V)}\,\forall\,\bw\,\in\,{\underline{\bb V}}\,,\\
{{(ii)}}& \quad \{ d_j \}_{j=1}^{\infty} \,\in\, 
l\hspace*{-.075cm}{l}^{{\bb W}}_{2 }\,\, \Leftrightarrow \, {\theta}\,\in\,
\W\,\,\,,
&&\;\|{\theta}\|_{D}\,=\,
\| \{ d_j \}_{j=1}^{\infty} \|_{2,(\bb W)}\,\forall\,\theta \,\in\,\W\,,
\end{align*}
 where we refer to the Dirichlet norms in (\ref{Dirichlet}). 
 \end{lem}
\noindent {\bf Proof.} Using our definitions, one easily sees that
\eqref{l2W_sequ} and \eqref{l2V_sequ} are the 
 properties of the Fourier coefficients for the energy spaces $\W$
 (for ${\boldsymbol L}$) and ${\underline{\bb V}}$ (for ${\boldsymbol
   S}$) (cf. \cite{Triebel} (4.1.8) and (4.1.9)). 
 \Qed
 \\[.2cm]
The Cartesian product space 
 $l\hspace*{-.075cm}{l}^{{\bb V}}_{2 }\,\times\,
l\hspace*{-.075cm}{l}^{{\bb W}}_{2 }$ 
is naturally equipped with the following scalar product and norm (cf.
Notation \ref{Cartesian_Pr}) 
\begin{align*} 
[(\mathbf{a},\mathbf{b}),\,(\mathbf{c},\mathbf{d})]_{l\hspace*{-.075cm}{l}^{{\bb V}}_{2 }\,\times\,
l\hspace*{-.075cm}{l}^{{\bb W}}_{2 }} 
&:= \,\big((\mathbf{a},\mathbf{c})_{2,(\bb V)}\,+\,(\mathbf{b},\mathbf{d})_{2,(\bb W)}\big)
\quad&
\forall\,(\mathbf{a},\mathbf{b}),\,(\mathbf{c},\mathbf{d})
\,\,\in\,l\hspace*{-.075cm}{l}^{{\bb V}}_{2 }\,\times\,
       l\hspace*{-.075cm}{l}^{{\bb W}}_{2 }\,\\
  \|({\mathbf{c}}\,.
{\mathbf{d}})\|_{l\hspace*{-.075cm}{l}^{{\bb V}}_{2 }\,\times\,
l\hspace*{-.075cm}{l}^{{\bb W}}_{2 }}&:=\,\sqrt{[(\mathbf{c},\mathbf{d}),
(\mathbf{c},\mathbf{d})]_{l\hspace*{-.075cm}{l}^{{\bb V}}_{2 }\,\times\,
l\hspace*{-.075cm}{l}^{{\bb W}}_{2 }}},.&
\end{align*}
After these preparations we are now in the position to carry out
the proof of Theorem \ref{T3.5}.\\[.2cm]
{\bf Proof.} {\em{(Theorem \ref{T3.5})}}\\
We rewrite Problem A in the weak formulation,
(where we use the potential $S$ from
\eqref{gradS}), i.e.
we are looking for 
$(\lambda,{\bw},{\theta})\in  
{\bb C}\times{\underline{\bb V}}\times {\W}
$ fulfilling for arbitrary ${\bu}\,\in\,{\underline{\bb V}}$ and
${\sigma}\,\in\,{\W}$ the equations
\begin{align}\label{bilinear_sepa1} 
\big({\theta},({\underline{\nabla}}\:\! S)^T{\bu}\big)
=\,\lambda\,\big({{\bu}},{\bw}\big)_D
\,\,,\\
\big({\sigma},({\underline{\nabla}}\:\! S)^T{\bw}\big)\,
=\,\lambda\,\big({\theta},
{\sigma}\big)_D\,\,.
\label{bilinear_sepa2} 
\end{align}
By taking the sum of \eqref{bilinear_sepa1} and  \eqref{bilinear_sepa2}
we obtain
\begin{align}\label{bilinear_weiter} 
\big({\sigma},({\underline{\nabla}}\:\! S)^T{\bw}\big)\,
+\,\big({\theta},({\underline{\nabla}}\:\! S)^T{\bu}\big)
=\,\lambda\Big(({{\bu}},{\bw})_D+({\theta},
{\sigma})_D\Big)\,\,,
\end{align}
which corresponds to to \eqref{bilinear_formulation} with $\lambda=\Ra_c$.
We express the elements ${\bu}$, ${\bw}$, ${\sigma}$
and ${\theta}$ as Fourier series in the 
systems of eigenfunctions $\{\bv_{j}\}_{j=1}^{\infty}$ and $\{\chi_{j}\}_{j=1}^{\infty}$ with respect to  of ${\underline{\bb H}}$ and ${\bb L}_{2}(.)$, respectively, i.e.
\begin{align}\label{biliFv}
{\bw}\,=\,\sum_{j=1}^{\infty} c_{j}\bv_{j}\,\quad\quad\,,\quad
{\bu}\,=\,\sum_{j=1}^{\infty} a_{j}\bv_{j}\,\,\\
{\theta}=\,\sum_{j=1}^{\infty}d_{j}\chi_{j}\,\quad\,,\quad
{\sigma}=\,\sum_{j=1}^{\infty}b_{j}\chi_{j}\quad .
\label{biliFT}
\end{align}
Moreover, we set
\begin{align}\label{Abbrsequ}
\mathbf{a}\,=\,\{ a_{j}\}_{j=1}^{\infty}\,,\,
\mathbf{b}\,=\,\{ b_{j}\}_{j=1}^{\infty}\,,\,
\mathbf{c}\,=\,\{ c_{j}\}_{j=1}^{\infty}\,\,\text{and}
\,\,
\mathbf{d}\,=\,\{ d_{j}\}_{j=1}^{\infty}
\end{align}
and regard them as sequences
of Fourier coefficients.
For sequences $\mathbf{a}$ considered as a row  we write in this
sense $\mathbf{a}^T$\,. In what follows we are going to use the usual
notation of double series as products  
of matrices. This method was employed, e.g.,~by E.~Schmidt \cite{Schmidt} and is a
standard method of functional analysis. The background for this method is 
the rearrangement theorem for double series. Inserting \eqref{biliFv}
and \eqref{biliFT} in the equations 
\eqref{bilinear_sepa1} and \eqref{bilinear_sepa2} results in
\begin{align}\label{BilinMatsequ1}
\mathbf{a}^T \big({C}\cdot
  \mathbf{d}\,-\,\lambda\,[\mbox{diag}(\kappa^2_{k})]\cdot\mathbf{c}\big) 
\,=\,0\,,\\
\label{BilinMatsequ2}
\mathbf{b}^T \big({C}^T\cdot
  \mathbf{c}\,-\,\lambda\,[\mbox{diag}(\omega^2_{k})]\cdot\mathbf{d}\big) 
\,=\,0\,\,,
\end{align}
where
\begin{equation}\label{symm1}
C_{j,k}:=\,(\chi_{k},({\underline{\nabla}}\:\! S)^T{\bv}_j)\,.
\end{equation}
One easily sees that $|C_{j,k}|\leq\, \sqrt{2}\,\cdot \gamma_{S}$ with
$\gamma_{S}\,:=\,\max_{\bx\in{\Omega_\mathcal{A}}}\big(|\frac{\partial
  S}{\partial r}|+|\frac{\partial S}{r\partial \varphi}|\big) 
$,
using that $\|{\underline{v}}_{j}\|_{{\underline{\bb 
H}}}\,=1$ and  $\|\chi_{k}\|_{{{\bb 
L}_{2}}}\,=1$ $\forall \,{j,k}\,\in\,{\bb N}\times{\bb N}$.
In $l\hspace*{-.075cm}{l}_{2}\times l\hspace*{-.075cm}{l}_{2}$ the sum
of the equations (\ref{BilinMatsequ1}) and (\ref{BilinMatsequ2}) 
can be written as
\begin{equation}\label{BLGLS2}
\big( \mathbf{a}^T, 
 \mathbf{b}^T
 \big)
\left[\left(\begin{array}{cc} 0&C\\
 C^T & 0
 \end{array}\right)-\lambda\,\left(\begin{array}{cc}\,\mbox{diag}(\kappa^2_{k}) &0\\
 0 &\,\mbox{diag}(\omega^2_{k})
 \end{array}\right)\right]\left(\begin{array}{c} \mathbf{c}\\
 \mathbf{d}
 \end{array}\right)=0\, ,
 \end{equation}
i.e. the sequences of Fourier coefficients corresponding to ${\bu}$, ${\bw}$, ${\sigma}$
and ${\theta}$ belong to 
$l\hspace*{-.075cm}{l}^{{\bb V}}_{2}$ and $l\hspace*{-.075cm}{l}^{{\bb W}}_{2}$.
Also, ${\tilde{\mathbf{a}}}$, ${\tilde{\mathbf{b}}}$, 
${\tilde{\mathbf{c}}}$, ${\tilde{\mathbf{d}}}\,\in\,l\hspace*{-.075cm}{l}_{2}$, with
\begin{align}\label{AbbrsequWeight} 
{\tilde{\mathbf{a}}}\,=\,\{ {\tilde{a}}_{j}:=\kappa_{j}a_{j}\}_{j=1}^{\infty}\,,\,
{\tilde{\mathbf{b}}}\,=\,\{ {\tilde{b}}_{j}:=\omega_{j}b_{j}\}_{j=1}^{\infty}\,,\,
{\tilde{\mathbf{c}}}\,=\,\{ {\tilde{c}}_{j}:=\kappa_{j}c_{j}\}_{j=1}^{\infty}\,\,\text{and}
\,\,
{\tilde{\mathbf{d}}}\,=\,\{ {\tilde{d}}_{j}:=\omega_{j}d_{j}\}_{j=1}^{\infty}\,,
\end{align}
because of ${\bu}, \,{\bw}\, \in\,{\underline{\bb V}}$ and  ${\sigma}
,\,{\theta}\, \in\,\W$.
Additionally, we define the matrix ${\tilde{C}}$ by 
\begin{equation}\label{symmvollst}
{\tilde{C}}_{j,k}:=\,\frac{1}{\kappa_{j}\omega_{k}}
(\chi_{k},({\underline{\nabla}}\:\! S)^T{\bv}_j)\,=\,\,\frac{1}{\kappa_{j}\omega_{k}}
C_{j,k}
\,\quad \forall \,{j,k}\,\in\,{\bb N}\times{\bb N}.
\end{equation}
So now 
we are able to rewrite problem (\ref{BLGLS2}) in the 
$l\hspace*{-.075cm}{l}_{2}\times l\hspace*{-.075cm}{l}_{2}$-setting, namely,
 \begin{equation}\label{BLGLS3}
\big( {\tilde{\mathbf{a}}}^T, 
 {\tilde{\mathbf{b}}}^T
 \big)
\left[\left(\begin{array}{cc} 0&{\tilde{C}}\\
 {\tilde{C}}^T & 0
 \end{array}\right)-\lambda\,\left(\begin{array}{cc}\,I &0\\
 0 &\,I
 \end{array}\right)\right]\left(\begin{array}{c} {\tilde{\mathbf{c}}}\\
 {\tilde{\mathbf{d}}}
 \end{array}\right)=0\, ,
 \end{equation} 
 where $I$ denotes the identity in  $l\hspace*{-.075cm}{l}_{2}$. 
We interpret (\ref{BLGLS3}) as an operator equation in 
$l\hspace*{-.075cm}{l}_{2}\times l\hspace*{-.075cm}{l}_{2}$. This is possible 
since the relation (\ref{BLGLS3}) is valid for arbitrary  $({\tilde{\mathbf{a}}},
{\tilde{\mathbf{b}}})\,\in l\hspace*{-.075cm}{l}_{2}\times
l\hspace*{-.075cm}{l}_{2}$. 
We define the linear and bounded operators 
$
{\tilde{{\cal C}}},\,{\cal I}: l\hspace*{-.075cm}{l}_{2}\times l\hspace*{-.075cm}{l}_{2}\longrightarrow  l\hspace*{-.075cm}{l}_{2}\times l\hspace*{-.075cm}{l}_{2}$ by
\begin{equation}\label{OPDEF}
{\tilde{{\cal C}}}\,:=\,
\left(\begin{array}{cc} 0&{\tilde{C}}\\
 {\tilde{C}}^T & 0
 \end{array}\right)\,\quad\,,\quad\,
{\cal I}\,:=\, 
\left(\begin{array}{cc}\,I &0\\
 0 &\,I
 \end{array}\right)\, .
 \end{equation} 
 Finally, equation (\ref{BLGLS3}) here is equivalent to the eigenvalue problem: 
 \begin{equation}\label{EVPlambda}
{\tilde{{\cal C}}}\left(\begin{array}{c} {\tilde{\mathbf{c}}}\\
 {\tilde{\mathbf{d}}}
 \end{array}\right)\,=\,
\lambda\,{\cal I}\left(\begin{array}{c} {\tilde{\mathbf{c}}}\\
 {\tilde{\mathbf{d}}}
 \end{array}\right)\,=\,\lambda\, \left(\begin{array}{c} {\tilde{\mathbf{c}}}\\
 {\tilde{\mathbf{d}}}
 \end{array}\right).
 \end{equation}
 The operators ${\tilde{{\cal C}}}$ and ${\cal I}$ are elements of the Banach space
 ${\cal L}\hspace*{-.19cm}{{\cal L}}(l\hspace*{-.075cm}{l}_{2}\times
 l\hspace*{-.075cm}{l}_{2},l\hspace*{-.075cm}{l}_{2}\times
 l\hspace*{-.075cm}{l}_{2})$ of linear and bounded operators.  We are
 now interested in the properties and the operator norms of
 ${\tilde{{\cal C}}}$ and  ${\cal I}$. It is well known, that the identity  ${\cal I}$ in $l\hspace*{-.075cm}{l}_{2}\times l\hspace*{-.075cm}{l}_{2}$ is a self-adjoint operator in  ${\cal L}\hspace*{-.19cm}{{\cal L}}(l\hspace*{-.075cm}{l}_{2}\times l\hspace*{-.075cm}{l}_{2},l\hspace*{-.075cm}{l}_{2}\times l\hspace*{-.075cm}{l}_{2})$ and 
 $\|{\cal I}\|_{{\cal L}\hspace*{-.14cm}{{\cal L}}}\,=\,1$.
 The boundedness of the operator ${\tilde{{\cal C}}}$ is obvious (cf. (\ref{symm1}) and
 (\ref{symmvollst})). 
 To show the compactness of ${\tilde{{\cal C}}}$ we use that it can be approximated by a sequence
 of finite operators $\{{\tilde{{\cal C}_{\ell}}}\}_{\ell \in {\bb N}}$, 
 $\{{\tilde{{\cal C}_{\ell}}}\}_{\ell \in {\bb N}}\,\in\,{\cal L}\hspace*{-.19cm}{{\cal L}}(l\hspace*{-.075cm}{l}_{2}\times l\hspace*{-.075cm}{l}_{2},l\hspace*{-.075cm}{l}_{2}\times l\hspace*{-.075cm}{l}_{2})$, with $\{{\tilde{{\cal C}_{\ell}}}\}_{\ell \in {\bb N}}\,\xrightarrow [{\cal L}\hspace*{-.14cm}{{\cal L}}]{} \,
 {\tilde{{\cal C}}}$, where $\forall \,\ell\,\in \,{\bb N}$ we use  the matrices
\begin{equation}\label{finite}
{\tilde{C}}_{{\ell},j,k}:=\,\left\{\begin{array}{cc}\quad
\,{\tilde{C}}_{j,k}\,\,\forall \,{j,k}\,\in\,{\bb N}\times{\bb N}:&
{j,k}\,\leq\,{\ell}\\
{~} & {~} \\ \quad\,
0 \,\quad\forall \,{j,k}\,\in\,{\bb N}\times{\bb N}:& \,j\,>\,{\ell}\,\vee\, \,k\,>\,{\ell}
\end{array}\right.
\end{equation}
and $\forall \,\ell\,\in \,{\bb N}$ the definitions:
\begin{equation}\label{OPFinApp}
{\tilde{{\cal C}}}_{\ell}\,:=\,
\left(\begin{array}{cc} 0&{\tilde{C}_{\ell}}\\
 {\tilde{C}_{\ell}}^T & 0
 \end{array}\right)\,.
 \end{equation}
 The convergence $\{{\tilde{{\cal C}_{\ell}}}\}_{\ell \in {\bb N}}\,\xrightarrow [{\cal L}\hspace*{-.14cm}{{\cal L}}]{} \, {\tilde{{\cal C}}}$ follows from \eqref{symm1}, \eqref{symmvollst},  $\lim_{j\rightarrow\infty} \kappa_{j}=\infty\,$, and $\lim_{k\rightarrow\infty}\omega_{k}=\infty\,$.\\
Since ${\tilde{{\cal C}}}$ is evidently symmetric,
 one has to show that it is self-adjoint.
 This follows since ${D}(\tilde{{\cal C}})=l\hspace*{-.075cm}{l}_{2}\times l\hspace*{-.075cm}{l}_{2}$
 for the symmetric operator $\tilde{{\cal C}}$ (cf. \cite[Th.1 in 
 4.1.6]{Triebel}) 
.\\
The existence of all real eigenvalues ${\lambda}$ (and all the eigen-elements
 $({\bw},{\theta})\in  {\underline{\bb V}}\times {\W}$) is a simple 
consequence of \cite[Th.1 in 4.2.6]{Triebel} with regard to the properties 
 of ${\tilde{{\cal C}}}$ and the problem (\ref{EVPlambda}). One has to
 note, that  ${\tilde{{\cal C}}}$ as a linear and compact operator on
 the separable infinite-dimensional space
 $l\hspace*{-.075cm}{l}_{2}\times l\hspace*{-.075cm}{l}_{2}$ has at
 most a countably 
infinite number of eigenvalues, where the non-zero eigenvalues $\{\lambda_{j}\}_{j=1}^{\infty}$ 
(counted
in multiplicity) can be ordered by their absolute values (cf. \cite{Triebel} (Th.1, 4.2.6) ).
There we have to take into account the finite multiplicity of the
non-zero eigenvalues. The number 
zero: $0\,\in\,{\bb C}$ is an element of the continuous	spectrum
in either case as the only accumulation point of the eigenvalues of ${\tilde{{\cal C}}}$,
but not an eigenvalue of ${\tilde{{\cal C}}}$, in what follows, that ${\tilde{{\cal C}}}^{-1}$
exists as a linear but unbounded operator. So we finish the proof by showing that $\lambda =0$ is not an 
eigenvalue of ${\tilde{{\cal C}}}$ in the following Lemma \ref{LambdaNicht0}.\,\Qed\\[-.3cm]
\begin{lem}\label{LambdaNicht0}
$\lambda =0$ is not an eigenvalue of ${\tilde{{\cal C}}}$, but as an accumulation point of eigenvalues the only element of the continuous spectrum of {~}\,${\tilde{{\cal C}}}\,\in\,
 {\cal L}\hspace*{-.19cm}{{\cal L}}(l\hspace*{-.075cm}{l}_{2}\times l\hspace*{-.075cm}{l}_{2},l\hspace*{-.075cm}{l}_{2}\times l\hspace*{-.075cm}{l}_{2})$.
\end{lem}
\noindent {\bf Proof.} Assuming, that $\lambda =0$ is an eigenvalue of ${\tilde{{\cal C}}}$, resp. of
problem (\ref{E1})-(\ref{E2}), than eq. (\ref{E1})  becomes 
\[\theta{\underline{\nabla}} S =-{\underline{\nabla}} p\,.\]
By applying  the curl operator to both sides we obtain
\[{\underline{\nabla}}{{\times}}(\theta{\underline{\nabla}} S)={\underline{\nabla}}\theta{{\times}}{\underline{\nabla}} S =\underline{0}\,.\]
But, as it is easy to check,  the condition that ${\underline{\nabla}}\theta$ and ${\underline{\nabla}} S$ are parallel
 vectors requires that
 \[\theta =\frac{1}{b}\ln r +r\sin\varphi\]
(modulo some constants), which does not obey the boundary conditions for the temperature.
Finally we get by application of \cite{Triebel} (Th.1, 4.2.6).
the property of $\lambda =0$ to be the only element of the continuous spectrum of {~}\,${\tilde{{\cal C}}}$.\,\Qed
%%%%%%%%%%%%%
\section{The numerical approximation of the critical constant}
\label{NumApp}
%%%%%%%%%%%%%%
\noindent Taking into account, that  
the maximum of our functional
corresponds to the largest positive eigenvalue $\lambda=\lambda_{max}\,=\,\frac{2}{{\Ra}_c}$
of
problem (\ref{E1})-(\ref{E2}), respectively to the operator norm of ${\tilde{{\cal C}}}$,
we are going to create an approximation method. 
We use determinants of square matrices of finite dimension.
(regarded as elements of the sequence
of finite operators) to create approximations for $\lambda$. Following
ideas and methods of E. Schmidt  
\cite{Schmidt} it is easy to prove, that our sequence (or subsequence) of approximations 
$\{\lambda^{{{\ell}_{\tilde{j}}}}_{\mathcal{A}}\}_{\tilde{j}=1}^{\infty}\,=\,
\{\lambda^{{{\ell}_{\tilde{j}}}}\}_{\tilde{j}=1}^{\infty}$ is 
monotonically increasing with the limit $\lambda\,=\,\|{\tilde{{\cal C}}}\|_{{\cal L}\hspace*{-.14cm}{{\cal L}}}$.\\
One could study the equation for the eigenvalues $\lambda$ of 
${\tilde{{\cal C}}}$ also in the sense of a limes of 	
finite-dimensional ${2\ell}\times{2\ell}$-matrices 
$\{\Gamma_{\ell}\}_{\ell \in {\bb N}}\,,\,\{{\tilde{{\cal J}_{\ell}}}\}_{\ell \in {\bb N}}
 \,$, 
 with ${\Gamma_{\ell}}\,{{\cal J}_{\ell}},\forall \,\ell \,\in {\bb N}\,$ according to description 
 \eqref{BLOMAinApp}:
 \begin{equation}\label{detA2}
 \lim_{\ell \to \infty} 
 \text{det}
 \left[{\tilde{{\Gamma}_{\ell}}} -\lambda\,{{{\cal J}_{\ell}}} \right]\,
 =:\text{det}
 \left[{\tilde{{\cal C}}} -\lambda\,{{{\cal I}}} \right]\,=\,
 \text{det}
 \left[\left(\begin{array}{cc} 0&{\tilde{C}}\\
 {\tilde{C}}^T & 0
 \end{array}\right)-\lambda\,\left(\begin{array}{cc}\,I &0\\
 0 &\,I
 \end{array}\right)\right]\, =\,0\,.
 \end{equation}
 {\bf Remark:} We choose a number $\ell\,\in\,{\bb N} $ in such a way,
 that the truncated eigenpair-systems of the Laplacian and the Stokes operator written as
\begin{equation} \label{Tranceigenpairs}
\,\,\{(\omega_{j})^2,\chi_{j}\}_{j=1}^{\ell}\, 
\,\qquad \text{and}\,\,\qquad
\{(\kappa_{j})^2,{\underline{v}}_{j}\}_{j=1}^{\ell}
\end{equation}
have the following properties: 
\begin{enumerate}
\item 
The spaces spanned by the systems $\{\chi_{j}\}_{j=1}^{\ell}$ and 
$\{{\underline{v}}_{j}\}_{j=1}^{\ell}$ contain the
entire eigenspaces of all eigenvalues $\{(\omega_{j})^2\}_{j=1}^{\ell}$ resp.
$\{(\kappa_{j})^2\}_{j=1}^{\ell}$
, particularly of $(\omega_{\ell})^2$ and 
$(\kappa_{\ell})^2$.
\item The number $\ell$ is specified by explicite determination of
the eigenpair-systems until a number $m\,\gg \ell$.\end{enumerate}
\begin{lem} \label{lsequence}
There exists  a	
strictly monotonically increasing sequence $\{{\ell}_{\tilde{j}}\}_{{\tilde{j}}=1}^{\infty}\,\longrightarrow \,\infty$
fulfilling with (\ref{Tranceigenpairs}) as:
\begin{equation} \label{SubTranceigenpairs}
\,\,\{(\omega_{j})^2,\chi_{j}\}_{j=1}^{\ell_{\tilde{j}}}\, 
\,\qquad \text{and}\,\,\qquad
\{(\kappa_{j})^2,{\underline{v}}_{j}\}_{j=1}^{\ell_{\tilde{j}}}
\end{equation}
the properties of the remark above for any ${\cal A}\,\in\,(0,\infty)$.
\end{lem}
{\bf Proof.} 
We use the intermeshing property of eigenvalues $\{(\omega_{j})^2\}_{j=1}^{\infty}$ resp.
$\{(\kappa_{j})^2\}_{j=1}^{\infty}$ which grow in the exact same
manner with respect to
their multiplicities. There one has to note, that we have ${\Omega}_{\cal A}$ as a 2d-domain 
and that the divergence $\di\bv_{j}\,=\,0\,$ works like a one-dimensional restriction for
2d-vector fields on 2d-domains.
\,\Qed\\[2mm]
For the approximations of 
$\lambda\,=\,\|{\tilde{{\cal C}}}\|_{{\cal L}\hspace*{-.14cm}{{\cal
      L}}}$ we use the Courant minimax principle. 
Let us introduce step by step notations for sets of square block matrices $\forall \,\ell\,\in \,\{{\ell}_{\tilde{j}}\}_{{\tilde{j}}=1}^{\infty}\,\subset\,{\bb N}$: 
We denote by  ${\hat{C}}_{\ell}$ the ${\ell}\times{\ell}$-matrices 
${\hat{C}}_{\ell}$ \,:
\begin{equation}\label{hatfinite}
{\hat{C}}_{{\ell},j,k}:=\,
\,{\tilde{C}}_{j,k}\,\,\forall \,{j,k}\,\in\,{\bb N}\times{\bb N}:\,
{j,k}\,\leq\,{\ell}\,\,,
\end{equation}
where ${\tilde{C}}_{j,k}$ is defined in (\ref{symmvollst}). 
The identity matrix of size ${\ell}$  is denoted by
${\hat{I}}_{\ell}$. ${\hat{I}}_{\ell}$ is the known square matrix with  
\begin{equation}\label{einhfinite}
{\hat{I}}_{{\ell},j,k}:=\,
\,{\delta}_{j,k}\,\,\forall \,{j,k}\,\in\,{\bb N}\times{\bb N}:\,
{j,k}\,\leq\,{\ell}\,\,,
\end{equation}
where ${\delta}_{j,k}$ is the Kronecker delta ($\forall \,\ell\,\in
\,{\bb N}$). 
For all $\ell\,\in \,\{{\ell}_{\tilde{j}}\}_{{\tilde{j}}=1}^{\infty}$
we define the block matrices $\Gamma_{\ell}$ and $J_{\ell}$ by
\begin{equation}\label{BLOMAinApp}
\Gamma_{\ell}\,:=\,
\left(\begin{array}{cc} 0&{\hat{C}_{\ell}}\\
 {\hat{C}_{\ell}}^T & 0
 \end{array}\right)\,\quad\,,\quad
 J_{\ell}\,:=\,\hat{I}_{2 \ell}\,=\,
\left(\begin{array}{cc} {\hat{I}_{\ell}} & 0\\
0 & {\hat{I}_{\ell}}
 \end{array}\right)\,.
 \end{equation}
One can also introduce a sequence of block matrices $\{\Gamma_{\ell}\}_{\ell\,\in \,\{{\ell}_{\tilde{j}}\}_{{\tilde{j}}=1}^{\infty},{\ell}_{\tilde{j}}<{\ell}_{\tilde{j}+1}}$
from a theoretical point of view. \\
It is worth to note, that the numerical methods and computer algebra systems for the evaluation
of all eigenvalues are restricted to small matrix sizes we do
not try to compute all eigenvalues of the matrices
$\{\Gamma_{\ell}\}_{\ell\,\in \,\{{\ell}_{\tilde{j}}\}_{{\tilde{j}}=1}^{\infty},{\ell}_{\tilde{j}}<{\ell}_{\tilde{j}+1}}$
to focus on approximations of 
$\lambda\,=\,\|{\tilde{{\cal C}}}\|_{{\cal L}\hspace*{-.14cm}{{\cal L}}}$.
The Courant minimax principle for $\lambda\,=\,\|{\tilde{{\cal C}}}\|_{{\cal L}\hspace*{-.14cm}{{\cal L}}}$ is standard
\begin{equation}\label{EWMaxC} 
\lambda\,=\,\|{\tilde{{\cal C}}}\|_{{\cal L}\hspace*{-.14cm}{{\cal L}}}\,:=\, 
\max \limits_{\|({\tilde{\mathbf{c}}},
 {\tilde{\mathbf{d}}})\|_{l\hspace*{-.075cm}{l}_{2}\times l\hspace*{-.075cm}{l}_{2}}=1
 }
 \big[
 ({\tilde{\mathbf{c}}},
 {\tilde{\mathbf{d}}})\,\,
{\tilde{{\cal C}}}\left(\begin{array}{c} {\tilde{\mathbf{c}}}\\
 {\tilde{\mathbf{d}}}
 \end{array}\right) \big]\,\,\,,
\end{equation}
or again with Courant minimax written as a simple limiting
    value like in the proof of \cite[Th. in Sec. 4.2.5]{Triebel}
\begin{equation}\label{EWMaxell} 
\lambda\,=\,
\|{\tilde{{\cal C}}}\|_{{\cal L}\hspace*{-.14cm}{{\cal L}}}\,=\, \lim \limits_{\tilde{j}\to\infty}
\max \limits_{\bz\in{\bb R}^{2{\ell}_{\tilde{j}}}:\|\bz\|=1
 }
 (\bz^T\,\cdot
{\tilde{\Gamma}}_{{{\ell}_{\tilde{j}}}} \cdot\bz\big)\,\,\,,
\end{equation}
where $\| . \|$ denotes the Euclidian norm 
of the  ${\bb R}^{2{\ell}_{\tilde{j}}}$ .\\
One can use the straightforward method of the evaluation of all of
${\tilde{\Gamma}}_{{{\ell}_{\tilde{j}}}}$ for
${\tilde{j}}\,=\,1,\,2,\,3 $ to calculate  
\begin{equation}\label{EWMaxelljot} 
\lambda^{{{\ell}_{\tilde{j}}}}\,:=\,
\max \limits_{\bz\in{\bb R}^{2{\ell}_{\tilde{j}}}:\|\bz\|=1
 }
 (\bz^T\,\cdot
{\tilde{\Gamma}}_{{{\ell}_{\tilde{j}}}} \cdot\bz\big)\quad\quad \forall \,{\tilde{j}}\,\in\,
{\bb N}
\end{equation}
in a first step. The start process sector is accurately described by the order $2{\ell}_{\tilde{j}}$ of the square matrices through $2{\ell}_{\tilde{j}}\,\approx \,40$ .\\
We get the inequalities: 
$\lambda^{{{\ell}_{1}}}\,\leq\,\lambda^{{{\ell}_{2}}} \,\leq\,\lambda^{{{\ell}_{3}}}\,\,\leq\,\dots $ using 
Lemma \ref{lsequence} and the definition (\ref{EWMaxelljot}) of  
$\lambda^{{{\ell}_{\tilde{j}}}}$, where we recognize the sequence 
$\{\lambda^{{{\ell}_{\tilde{j}}}}\}_{\tilde{j}=1}^{\infty}$ to be monotonically increasing
with the limit $\lambda$ as a simple consequence. So we  approach	
the eigenvalue $\lambda$ from below. 
We take the square (block) matrices of order $2{\ell}_{\tilde{j}}$: 
$\Gamma_{{\ell}_{\tilde{j}}}\,$ and $\,J_{{\ell}_{\tilde{j}}}$ (cf. (\ref{BLOMAinApp}))
to explain the functions 
\begin{equation}\label{Functdet}
F_{\tilde{j}}(\lambda)\,:=\,
 \text{det}
 \left[\Gamma_{{\ell}_{\tilde{j}}}- \lambda\,J_{{\ell}_{\tilde{j}}} \right]\,\quad\,
 \forall\,{\tilde{j}}\,\in\,{\bb N}\,.
 \end{equation}	
We evaluate $\lambda^{{{\ell}_{\tilde{j}}}}$ by calculating the zeros of $F_{\tilde{j}}(\lambda)$ in the interval 
$[\lambda^{{{\ell}_{\tilde{j}-1}}}
,
\lambda^{{{\ell}_{\tilde{j}-1}}} +\delta(\cal A)]$. 
There we use the bisection method.
The initial approximations (start interval) are chosen by means of 
$\lambda^{{{\ell}_{\tilde{j}-1}}}$ to be
$[\lambda^{{{\ell}_{\tilde{j}-1}}}
,
\lambda^{{{\ell}_{\tilde{j}-1}}} +\delta(\cal A)]$, where one can take the constant 
$\delta({\cal A}) \,\approx \,    1$. The initial approximations are taken like 
in the complete induction method with ${\tilde{j}-1}\,\geq\,1$. 
This method for the approach of  $\lambda$ is also restricted by the matrix size.
Because of the very long calculation time for determinants for large
${{\ell}_{\tilde{j}}}$ we make a rough estimate for the border area
around  ${\ell}_{\tilde{j}}\,\approx \,400$. 
That applies to other representations of $F_{\tilde{j}}(\lambda)$ (cf.(\ref{Functdet})) like 
\begin{equation}\label{BlockFunctdet}
F_{\tilde{j}}(\lambda)\,:=\,
 \text{det}
 \left[{\lambda}^2\cdot{\hat{I}}_{{\ell}_{\tilde{j}}}  - {\hat{C}}_{\tilde{j}}\cdot {\hat{C}}_{\tilde{j}}^T \right]\,\quad\,
 \forall\,{\tilde{j}}\,\in\,{\bb N}\,.
\end{equation}	
\noindent Now we  present some results for the numerical approximation of critical Rayleigh numbers. The calculations were implemented and performed in Maple 2021 (and Maple 2022)
to ensure tight tolerances for the eigenfunctions, eigenvalues and
${{\Gamma}}_{{{\ell}_{\tilde{j}}}}$. 
In particular we use some tools for the investigation of
Besselfunctions (cf., e.g., \cite{Rummler} and 
\cite{RRTZAMM}).\\
We present our 
numerical results of approximations for ${\Ra}_c(\mathcal{A})$ in the cases $\mathcal{A} = 1$
and $\mathcal{A} = 10$ (cf. \cite{Markert}). 
The numbers $\ell_{j}$ are choosen in the way described
in Lemma \ref{lsequence}.\\[.2cm]
\begin{tabular}{|c|r|r|}
\hline 
$\ell_{j}$  & $\lambda_{1}^{\ell_{j}}$ \hspace*{4cm}&  ${\Ra}_c(\mathcal{A} = 1) \hspace*{4cm}$\\
\hline
7 &
 0.065275715709889022434641241556310&
 30.639265739938989496469183142650 \\
\hline
12&
 0.067433400164287488107980906050821&
 29.658892998535072667362035625305\\
\hline
22&
 0.070024363834146514766687200657022&
 28.561487609327265913546978946005\\
\hline
49&
 0.070752337753538312200877678082614&
 28.267617205340757028944998021332\\
\hline
92&
 0.070843842842887568707945584675261&
 28.231105481325421889287981461511\\
 \hline
300&
 0.070926398413629912903178212704811&
 28.198245571928834452591903765447\\
\hline 
\end{tabular}
{~}\\[.9cm]
\begin{tabular}{|c|r|r|}
\hline 
$\ell_{j}$  & $\lambda_{1}^{\ell_{j}}$ \hspace*{4cm}&  ${\Ra}_c(\mathcal{A} = 10) $\hspace*{4cm}\\
\hline
11&
 0.036922959760029763204380893713543&
 54.166838547029519643082315254994\\
\hline
21&
 0.036922959760029763204380893713509&
 54.166838547029519643082315255044\\
\hline
39&
 0.036988540331501778832593839394502&
 54.070800904156620330758041435409\\
\hline
62&
 0.036988763473199326352547015240593&
 54.070474711843912995766624606994\\
\hline
90&
 0.036988784077179802197566851096866&
 54.070444592794772101553849052213\\
 \hline
\end{tabular}
{~}\\
{~}\\[.3cm]
Here one can observere, that
the sequences 
$\{\lambda_1^{{{\ell}_{\tilde{j}}}}\}_{\tilde{j}\geq 1}$ are  monotonically increasing. This increasing behavior of the 
$\lambda_1^{{\ell}_{\tilde{j}}}(\mathcal{A})$ is
justified by the growing dimension of the  variation-sets (the matrix sizes of the 
${{\Gamma}}_{{{\ell}_{\tilde{j}}}}$) in both tables.
%%%%%%%%%%
\section{Alternative proof via the stream-function}
%%%%%%%%%%%%%%%%%%%%%%%%%%%%%%%%%%%%
We will sketch another way for a proof of a modified version of Theorem \ref{T3.5} using streamfunctions. 
For this reason we rephrase the 
propositions of Theorem \ref{T3.5} in the following Theorem \ref{SFTheor}.
\begin{theor}\label{SFTheor}
For all $\mathcal{A}$ Problem A (cf.~\eqref{E1}-\eqref{E2}) is 
equivalent to the eigenvalue 
problem of the compact self-adjoint operator\, ${\tilde{{\cal B}}}\,\in\,
 {\cal L}\hspace*{-.19cm}{{\cal L}}(l\hspace*{-.075cm}{l}_{2}\times l\hspace*{-.075cm}{l}_{2},l\hspace*{-.075cm}{l}_{2}\times l\hspace*{-.075cm}{l}_{2})$ (see below). 
All the eigenvalues $\lambda$
are real
and $\lambda\,=\,0$ is not an eigenvalue of the operator
$\widetilde{\cal B}$, but is as  an accumulation point of eigenvalues the only element of the continuous spectrum of {~}\,${\tilde{{\cal B}}}$.
The critical Rayleigh number ${\Ra}_c$ is fixed through
\[\qquad %\lambda\,=
  \frac{2}{{\Ra}_c}\,:=\, |\lambda_{max}|\,=\,\|{\tilde{{\cal B}}}\|_{{\cal L  }\hspace*{-.15cm}{{\cal L}}}\,\,.\]
The set $\{\lambda\}$ of the ${\tilde{{\cal B}}}$-eigenvalues
consists of an countably infinite number of (finite multiplicity) eigenvalues.\\ 
The eigen-triples
$(\lambda,{\tilde{\mathbf{a}}},
{\tilde{\mathbf{b}}})\,\in {\bb C}\times l\hspace*{-.075cm}{l}_{2}\times l\hspace*{-.075cm}{l}_{2}$
correspond to  a solution 
$(\lambda,{\tilde \bw},{\tilde \theta})\in  
{\bb C}\times{\underline{\bb V}}\times {\W}
$
of Problem A (cf.~\eqref{E1}-\eqref{E2}).
\end{theor}
\noindent {\bf Remark:} The use of streamfunctions will need additional
efforts in the treatment of 
Problem A. Especially we get worse estimates \eqref{eststream} in comparison to 
the proof of Theorem \ref{T3.5} %and an deteriorate convergence of the sequence of finite opererators to $\widetilde{\cal B}$ (at showing of the compactness of $\widetilde{\cal B}$) also.
 {~}\\[.2cm]
In a first step we are going to make some arrangements for the proof
of Theorem \ref{SFTheor}. Again, we use
 Fourier representations for any streamfunction
 $\Psi\,\in \,{{\bb W}_{2}^{2}\hspace{-.62cm}{~}^{{~}^{{~}^{o}}}}\,\,\,\,$ as
expansions in the systems of eigenfunctions
$\{\psi_{j}\}_{j=1}^{\infty}$ for the Bilaplacian ${\boldsymbol B}$. 
\begin{lem}
\label{CompletenessStr}
The system  (cf. Notation \ref{eigenpairs})
$
\left\{(\psi_{j},\chi_k)\right\}_{j,k \in \NN}\quad 
$
is a complete basis for the space
 ${{\bb
  W}_{2}^{2}\hspace{-.62cm}{~}^{{~}^{{~}^{o}}}}\,\,\,\times \,{{\bb
  W}_{2}^{1}\hspace{-.62cm}{~}^{{~}^{{~}^{o}}}}$\,\, .
\end{lem}
\noindent {\bf Proof.} We have to apply the properties of positive self-adjoint operators with a pure point
spectrum (cf. Theorem \ref{thmstok} and Lemma \ref{STOeiFU}, or \cite[Th. in 4.5.1]{Triebel}).
\Qed \\[.1cm]
%
%%%
In what follows we are going to use the Fourier representations
for any streamfunction $\Psi\,\in \,{{\bb W}_{2}^{2}\hspace{-.62cm}{~}^{{~}^{{~}^{o}}}}$ \;
and any temperature ${\theta}\,\in \,\W$
\, 
as expansions in the systems of eigenfunctions $\{\psi_{j}\}_{j=1}^{\infty}$ 
and $\{\chi_{j}\}_{j=1}^{\infty}$ in the sense of ${\bb L}_{2}$ respectively (cf. Lemma \ref{STOeiFU} and \eqref{eigenpairs2}):
\begin{align} \label{F-expensions1}
\Psi =\,\sum_{j=1}^{\infty} c_{j}\psi_{j}\,\qquad\qquad {\text{and}}\qquad\qquad {\theta}=\,\sum_{j=1}^{\infty}d_{j}\chi_{j}\,.
\end{align}
The sequences $\{c_{j}\}_{j=1}^{\infty}$ and $\{d_{j}\}_{j=1}^{\infty}$ are - as
sequences of real or complex numbers - elements of {\em Hilbert space of sequences}
$l\hspace*{-.075cm}{l}_{2}$ (\cite[Lemma 2. in 2.1.3]{Triebel}). Parseval's equation is fulfilled because of the properties of the complete orthonormal 
systems $\{\psi_{j}\}_{j=1}^{\infty}$ and $\{\chi_{j}\}_{j=1}^{\infty}$ (in ${\bb L}_{2}$).\,
(We note that the system $\{\psi_{j}\}_{j=1}^{\infty}$ is an orthogonal system in the space
${{\bb W}_{2}^{2}\hspace{-.62cm}{~}^{{~}^{{~}^{o^*}}}}\,\,:=\,
({{
  W}_{2}^{2}\hspace{-.62cm}{~}^{{~}^{{~}^{o}}}}\,\,\,,\|
. \|_{\Delta})$. Here the vector space $ {{
  W}_{2}^{2}\hspace{-.62cm}{~}^{{~}^{{~}^{o}}}}$ \,\,\,
is endowed with the norm to the ${\Delta}$-(quasi-)scalar product (cf. Notation \ref{W2Aequiv}). 
The system $\{\chi_{j}\}_{j=1}^{\infty}$ is a complete 
orthogonal system in  $\W$ by construction.
\begin{notation} \label{W2Aequiv} 
$\forall \,\Psi\,,\Phi \,\in  \, {{\bb
  W}_{2}^{2}\hspace{-.62cm}{~}^{{~}^{{~}^{o}}}}\,\,\,\,$ 
by
\begin{align}\label{Delta_Skal} 
({{\Phi}},{\Psi})_{\Delta}\,:=\,(\Delta{{\Phi}},\Delta{\Psi})
\end{align}
we denote the ${\Delta}$-(quasi-)scalar product on the vector space $ {{
  W}_{2}^{2}\hspace{-.62cm}{~}^{{~}^{{~}^{o}}}}$\,\,\,\,.\\
\end{notation}
{~}\\[-.5cm]
{\bf Remark:} The ${\Delta}$-(quasi-)scalar product generates an equivalent to the
standard norm of ${{\bb
  W}_{2}^{2}\hspace{-.62cm}{~}^{{~}^{{~}^{~}}}}\,\,\,\,$ 
on the vector space $ {{
  W}_{2}^{2}\hspace{-.62cm}{~}^{{~}^{{~}^{o}}}}$\,\,\,.
\begin{notation} \label{SequSpacesW2} 
{\em (Special sequence space of coefficients)}\\ 
Let us regard the sequence of eigenvalues from the  eigenpairs-system
for the Bilaplacian ${\boldsymbol B}$ 
(cf. Notation \ref{eigenpairs})
ordered by increasing eigenvalues taking into account their % \{\mu_{j}^2,\psi_{j}\}
multiplicities (and the correspondent eigenfunctions):
\begin{align} \label{eigenvalstream_sequ}
\{(\mu_{j})^2\,=\,|\mu_{j}|^2\}_{j=1}^{\infty}\,\quad \text{
and the correspondent sequences of positive roots:}
\quad
\{\mu_{j}\}_{j=1}^{\infty}\,\,.
\end{align}
We define the space $l\hspace*{-.075cm}{l}^{{\bb W}2}_{2 }$\,,\,
with 
respect to spectral operators in the following way: 
$l\hspace*{-.075cm}{l}^{{\bb W}2}_{2 }\,:=\,
(l^{{\bb W}2}_{2 },\| . \|_{2,(\bb W)})$ is defined by 
\begin{eqnarray}\label{l2W2_sequ}
l^{{\bb W}2}_{2 }\,:=\,
\{ \,\{ a_j \}_{j=1}^{\infty}\,\in\,l_{2 }:
\,\sum_{j=1}^{\infty}\mu_{j}^2|a_j|^{2}\,<\,\infty\} \hspace*{6cm} \nonumber \\
\,\text{with the norm and the scalar product}\, \,\forall\,
\{ a_j \}_{j=1}^{\infty} \,,\{ b_j \}_{j=1}^{\infty} \,\in\, 
l^{{\bb W}2}_{2 } \,:\hspace*{2cm}\\
\| \{ a_j \}_{j=1}^{\infty} \|_{2,({\bb W} 2)}\,\,:=\,
\sqrt{(\{ a_j \}_{j=1}^{\infty},\{ a_j \}_{j=1}^{\infty})_{2,({\bb W} 2)}}\,\quad, \quad\,
(\{ a_j \}_{j=1}^{\infty},\{ b_j \}_{j=1}^{\infty})_{2,({\bb W} 2)}\,:=\,
\sum_{j=1}^{\infty}\, \mu_{j}^2 a_j {\DarkGr{\overline{\Black b}}}_j   \,.
\nonumber
\end{eqnarray}
\end{notation} 
\noindent Like in Section \ref{Secvelocity} Lemma \ref{Hilbert-sequences} we
obtain the following
\begin{lem}
\label{HilbertStream-sequences}
Let the Fourier coefficients of $\Psi\,\in \,{{\bb W}_{2}^{2}\hspace{-.62cm}{~}^{{~}^{{~}^{o}}}}\,\,\,\,$with respect to the ${\bb L}_{2}$
orthonormal system $\{\psi_{j}\}_{j=1}^{\infty}$ given by 
$\{ c_j \}_{j=1}^{\infty} $.\,
Then the following statements are true:
\begin{align*}
{{(iii)}}& \quad \{ c_j \}_{j=1}^{\infty} \,\in\, 
l\hspace*{-.075cm}{l}^{{\bb W}2}_{2 }\,\, \Leftrightarrow \,
\Psi\,\in \,{{\bb W}_{2}^{2}\hspace{-.62cm}{~}^{{~}^{{~}^{o}}}}\,\,\,\quad \text{and}
& \sqrt{(\Psi,\Psi)_{\Delta}}\,=\,
\| \{ c_j \}_{j=1}^{\infty} \|_{2,({\bb W}2)}\,\forall\,\Psi\,\in
                                                                                        \,{{\bb W}_{2}^{2}\hspace{-.62cm}{~}^{{~}^{{~}^{o}}}}\quad ,
\end{align*}
 where we refer to the scalar product in \eqref{Delta_Skal}. 
 \end{lem}
\noindent {\bf Proof.} One only has to rewrite our definitions.\Qed\\[.1cm]
Let us now consider the representation of two-dimensional (or
three-dimensional) vector fields by means of the streamfunction
$\Psi=\Psi(r,\varphi)$, 
or, by considering the reference frame of toroidal vector fields here chosen:
\[{\bw}={\underline{\nabla}}\,\Psi {{\times}} \beeis\,,\]
where $\beeis$ denotes the unit vector in  the horizontal direction.
We are going to employ the common vector product of the three-dimensional space
to avoid additional notations as the wedge-product. So we understand our terms
in the sense of three-dimensional expression and our results also as
three-dimensional vector fields with nonvanishing  coefficients of $\beeis$  respectively
in the following passage.
One can calculate the curl of all terms in \eqref{E1} and \eqref{E2} by assuming higher regularity 
\begin{align*}
{\underline{\nabla}}\times ({\theta}{\underline{\nabla}} S)&= {\underline{\nabla}}
{\theta}{{\times}} ({\underline{\nabla}} S)
\\
{\underline{\nabla}}\times
(\lambda\Delta{\bw})
&=-\lambda\Delta^2 \Psi
\beeis\\ ({\underline{\nabla}}\times{\underline{\nabla}})\,{p}&=\underline{0}\,.
\end{align*}
Hence, Problem A provides 
\begin{align}\label{3.21}
\beeis^T\cdot ({\underline{\nabla}}{\theta}{{\times}}{\underline{\nabla}}
S) -\lambda\Delta^2\Psi  &=0\\ \beeis^T\cdot (
{\underline{\nabla}}\Psi {{\times}}{\underline{\nabla}} S) -\lambda\Delta{\theta}&=0\,.\label{3.22}
\end{align} %ZZZ
%Letzte Beziehung Spatprodukt also Determinanten-Rechnen
We are going to investigate \eqref{3.21} and \eqref{3.22} as bilinear forms again. 
So we seek for 
$(\lambda,{\bw}={\underline{\nabla}}\,\Psi {{\times}} \beeis,{\theta})\in  
{\bb C}\times{\underline{\bb V}}\times {\W}
$ fulfilling for arbitrary ${\bu}={\underline{\nabla}}\,\Phi {{\times}} \beeis,$ and ${\sigma}$,
the equations (written with 
$\Psi\,,\Phi \,\in  \, {{\bb
  W}_{2}^{2}\hspace{-.62cm}{~}^{{~}^{{~}^{o}}}}\,\,\,\,$):
\begin{align}\label{bilstream_sepa1} 
\big({\Phi},\beeis^T\cdot ({\underline{\nabla}}{\theta}{{\times}}{\underline{\nabla}}
\:\! S )\big)
=\,\lambda\,\big({{\Phi}},{\Psi}\big)_{\Delta}
\,\,\,,\\\
\big({\sigma},\beeis^T\cdot ({\underline{\nabla}}{\Psi}{{\times}}{\underline{\nabla}}
\:\! S )\big)\,
=\,- \lambda\,\big(
{\sigma},{\theta}\big)_D\,\,,
\label{bilstream_sepa2} 
\end{align}
where we have used ${\bw}={\underline{\nabla}}\,\Psi {{\times}} \beeis \in {\underline{\bb V}}$ (cf. Problem A)
and  the ${\Delta}$-(quasi-)scalar product in \eqref{bilstream_sepa1}.
\begin{lem}\label{Schiefsymm} $\forall \,\psi\,,\chi \,\in  \,D({\boldsymbol L}) $, 
($\text{resp.}
\,\,\forall \,\psi\,,\chi \in  \,\W$) we have
\begin{equation}\label{3.23}
\big({\psi},\beeis^T\cdot ({\underline{\nabla}}\:\!{\chi}{{\times}}{\underline{\nabla}}
\:\! S )\big)\,=\,
\int_{\Omega_\mathcal{A}} (\beeis^T\cdot
{\underline{\nabla}}\:\!{\chi}{{\times}}{\underline{\nabla}}\:\!
S){\psi}d\:\!\bx=-\int_{\Omega_\mathcal{A}} (\beeis^T\cdot
{\underline{\nabla}}\:\!{\psi}{{\times}}{\underline{\nabla}} \:\!S)\chi d\:\!\bx\, 
\,=\,-\big({\chi},\beeis^T\cdot ({\underline{\nabla}}\:\!{\psi}{{\times}}{\underline{\nabla}}
\:\! S )\big).
\end{equation}
\end{lem}
\noindent {\bf Proof.}
Since ${\underline{\nabla}}\times{\underline{\nabla}} \:\!S = {\underline{0}}$, the left-hand side is 
equal to
\[\int_{\Omega_\mathcal{A}} \beeis^T\cdot
({\psi}{\underline{\nabla}}\times(\chi{\underline{\nabla}} S))d\:\!\bx=
\int_{\Omega_\mathcal{A}} \beeis^T\cdot
{\underline{\nabla}}\times(\psi\chi{\underline{\nabla}} S)d\:\!\bx -
\int_{\Omega_\mathcal{A}} \beeis^T\cdot
({\underline{\nabla}}\psi{{\times}}(\chi{\underline{\nabla}} S))d\:\!\bx
\]
and, on the other hand, here the first integral on the right-hand
side vanishes, because of the Gauss theorem and of the boundary
condition for $\psi$:
\[\int_{\Omega_\mathcal{A}}
\beeis^T\cdot {\underline{\nabla}}\times(\psi\chi{\underline{\nabla}} S)d\:\!\bx
=\int_{\partial \Omega_\mathcal{A}} \psi\chi ({\underline{\nabla}}\:\! S)^T
 \:\! d\bx\, =\,0.
\] 
{~}\\[-.85cm]
\Qed \\
 {~}\\[.25cm]
We are now in the position to give a shortened 
proof of Theorem \ref{SFTheor}.\\[.1cm]
{\bf Proof.} {\em{(Theorem \ref{SFTheor})}}\\
We start with a weak formulation
(cf. \eqref{E1}-\eqref{E2}) with the 
streamfunction.
By taking the sum of \eqref{bilstream_sepa1} and  \eqref{bilstream_sepa2}
we get:
\begin{align}\label{bilstream_weiter}
\big({\Phi},\beeis^T\cdot ({\underline{\nabla}}{\theta}{{\times}}{\underline{\nabla}}
\:\! S )\big) \,+\, \big({\sigma},\beeis^T\cdot ({\underline{\nabla}}{\Psi}{{\times}}{\underline{\nabla}}
\:\! S )\big)
=\,\lambda\,\Big(\big({{\Phi}},{\Psi}\big)_{\Delta}
\,-\,({\theta},
{\sigma})_D\Big)\,.
\end{align}
Now we consider the elements ${\Phi}$, ${\Psi} \,\in \, {{
  W}_{2}^{2}\hspace{-.62cm}{~}^{{~}^{{~}^{o}}}}\,\,\,\,$ and
${\sigma}$, ${\theta}\,\in \, \W\,$ as Fourier series in the 
systems of eigenfunctions $\{\psi_{j}\}_{j=1}^{\infty}$ and $\{\chi_{j}\}_{j=1}^{\infty}$ in the sense of ${\bb L}_{2}$ and ${\bb L}_{2}$ respectively. We set
\begin{align}\label{biliFSt}
{\Psi} \,=\,\sum_{j=1}^{\infty} c_{j}\psi_{j}\,\quad\quad\,,\quad
{\Phi}\,=\,\sum_{j=1}^{\infty} a_{j}\psi_{j}\,\,\\
{\theta}=\,\sum_{j=1}^{\infty}d_{j}\chi_{j}\,\quad\quad\,,\quad
{\sigma}=\,\sum_{j=1}^{\infty}b_{j}\chi_{j}\quad .
\label{biliFTS}
\end{align}
The sequences of Fourier coefficients are denoted by $\mathbf{a}$, $\mathbf{b}$, $\mathbf{c}$
and $\mathbf{d}$ as in 
\eqref{Abbrsequ}.
The sequences are regarded as 
columns. 
The Fourier series are
substituted in \eqref{bilstream_sepa1} and \eqref{bilstream_sepa2} and
\begin{align}\label{A}
\mathbf{a}^T \big({B} \mathbf{d}\,-\,\lambda\,[\mbox{diag}(\mu^2_{k})]\mathbf{c}\big)
\,=\,0\,,\\
\label{B}
\mathbf{b}^T \big({B}^T\mathbf{c}\,-\,\lambda\,[\mbox{diag}(\omega^2_{k})]\mathbf{d}\big)
\,=\,0\,\,,
\end{align}
where
\begin{equation}\label{Streamsymm1}
B_{j,k}:=\,({\psi}_j,\beeis^T\cdot ({\underline{\nabla}}\:\!{\chi}_k{{\times}}{\underline{\nabla}}
\:\! S ))\,
\,\quad \,,\, \,\,\quad \forall \,{j,k}\,\in\,{\bb N}\times{\bb N}.
\end{equation} 
The sum of the equations \eqref{A} and \eqref{B} (the equivalent to \eqref{bilstream_weiter})
can be written in $l\hspace*{-.075cm}{l}_{2}\times l\hspace*{-.075cm}{l}_{2}$ as:
\begin{equation}\label{streamLGLS2}
\big( \mathbf{a}^T, 
 \mathbf{b}^T
 \big)
\left[\left(\begin{array}{cc} 0&B\\
 B^T & 0
 \end{array}\right)-\lambda\,\left(\begin{array}{cc}\,\mbox{diag}(\mu^2_{k}) &0\\
 0 &\,\mbox{diag}(\omega^2_{k})
 \end{array}\right)\right]\left(\begin{array}{c} \mathbf{c}\\
 \mathbf{d}
 \end{array}\right)=0\, .
 \end{equation}
For Lemma \ref{Schiefsymm} we get the structure and the symmetry of the following
linear
operator written as a ``block matrix"{} on $l\hspace*{-.075cm}{l}_{2}\times l\hspace*{-.075cm}{l}_{2}$
\begin{equation*}
\left(\begin{array}{cc} 0&B\\
 {B}^T & 0
 \end{array}\right)\, .
\end{equation*}
But here  the elements $B_{j,k}$ in \eqref{Streamsymm1} are not
uniformly bounded in contrast to the uniform boundedness of 
the elements $C_{j,k}$ of $C$ in
Section \ref{Secvelocity} in
\eqref{symm1}.
In order to estimate the $B_{j,k}$ we need the equivalent norms on $ {{
  W}_{2}^{2}\hspace{-.62cm}{~}^{{~}^{{~}^{o}}}}\,\,\,\,$ respectively on $\W$. 
$\forall \,{j,k}\,\in\,{\bb N}\times{\bb N}$ this results in
\begin{equation}\label{eststream}
|B_{j,k}|=|({\psi}_j,\beeis^T\cdot ({\underline{\nabla}}\:\!{\chi}_k{{\times}}{\underline{\nabla}}
\:\! S ))|\,=\,|({\chi}_k,\beeis^T\cdot ({\underline{\nabla}}\:\!{\psi}_j{{\times}}{\underline{\nabla}}
\:\! S ))|\leq\, \gamma_{*}\gamma^{II}_{S}\,\,\max \big(\mu_{j},\omega_{k}\big)\,,
\end{equation}
where we have used $\gamma_{*}$ as constant for the two equivalent norms and
$\gamma^{II}_{S}\,:=\,\max_{\bx\in{\Omega_\mathcal{A}}}\sqrt{
\big(|\frac{\partial S}{\partial r}|^2+|\frac{\partial S}{r\partial \varphi}|^2\big)}
$,
$\|{\psi}_{j}\|_{{{\bb 
L}_{2}}}\,=1$, and  $\|\chi_{k}\|_{{{\bb 
L}_{2}}}\,=1$ $\forall \,{j,k}\,\in\,{\bb N}\times{\bb N}$.\\[.1cm]
Now the sequences of Fourier coefficients corresponding to  ${\Phi}$, ${\Psi}$
, ${\sigma}$
and ${\theta}$ belong to the spaces 
$l\hspace*{-.075cm}{l}^{{\bb W}2}_2$ and $l\hspace*{-.075cm}{l}^{{\bb W}}_{2}$
(cf. Lemma \ref{Hilbert-sequences} and Notation \ref{SequSpacesW2}).
Thus, it follows that ${\tilde{\mathbf{a}}}$, ${\tilde{\mathbf{b}}}$, 
${\tilde{\mathbf{c}}}$, ${\tilde{\mathbf{d}}}\,\in\,l\hspace*{-.075cm}{l}_{2}$, with
\begin{align}\label{AbbrsequWstream} 
{\tilde{\mathbf{a}}}\,=\,\{ {\tilde{a}}_{j}:=\mu_{j}a_{j}\}_{j=1}^{\infty}\,,\,
{\tilde{\mathbf{b}}}\,=\,\{ {\tilde{b}}_{j}:=\omega_{j}b_{j}\}_{j=1}^{\infty}\,,\,
{\tilde{\mathbf{c}}}\,=\,\{ {\tilde{c}}_{j}:=\mu_{j}c_{j}\}_{j=1}^{\infty}\,\,
\text{and}
\,\,
{\tilde{\mathbf{d}}}\,=\,\{ {\tilde{d}}_{j}:=\omega_{j}d_{j}\}_{j=1}^{\infty}\,.
\end{align}
Additionally, we define the matrix ${\tilde{B}}\,$ by 
\begin{equation}\label{streamvollst}
{\tilde{B}}_{j,k}:=\,\,\frac{1}{\mu_{j}\omega_{k}}
B_{j,k}
\,\quad \forall \,{j,k}\,\in\,{\bb N}\times{\bb N}.
\end{equation}
So now we are able to formulate the problem \eqref{streamLGLS2} equivalently in the 
$l\hspace*{-.075cm}{l}_{2}\times l\hspace*{-.075cm}{l}_{2}$-sense through
 \begin{equation}\label{BLGSF3}
\big( {\tilde{\mathbf{a}}}^T, 
 {\tilde{\mathbf{b}}}^T
 \big)
\left[\left(\begin{array}{cc} 0&{\tilde{B}}\\
 {\tilde{B}}^T & 0
 \end{array}\right)-\lambda\,\left(\begin{array}{cc}\,I &0\\
 0 &\,I
 \end{array}\right)\right]\left(\begin{array}{c} {\tilde{\mathbf{c}}}\\
 {\tilde{\mathbf{d}}}
 \end{array}\right)=0\, ,
 \end{equation} 
 where $I$ denotes the identity in  $l\hspace*{-.075cm}{l}_{2}$ again. 
We interpret (\ref{BLGSF3}) as an operator equation in 
$l\hspace*{-.075cm}{l}_{2}\times l\hspace*{-.075cm}{l}_{2}$. This is possible 
since the relation (\ref{BLGSF3}) is valid for arbitrary  $({\tilde{\mathbf{a}}},
{\tilde{\mathbf{b}}})\,\in l\hspace*{-.075cm}{l}_{2}\times l\hspace*{-.075cm}{l}_{2}$.\\
We define the linear and bounded operators ${\tilde{{\cal B}}}$ and ${\cal I}$ by;\quad 
$
{\tilde{{\cal B}}},\,{\cal I}: l\hspace*{-.075cm}{l}_{2}\times l\hspace*{-.075cm}{l}_{2}\longrightarrow  l\hspace*{-.075cm}{l}_{2}\times l\hspace*{-.075cm}{l}_{2}$ by:
\begin{equation}\label{SFOPDEF}
{\tilde{{\cal B}}}\,:=\,
\left(\begin{array}{cc} 0&{\tilde{B}}\\
 {\tilde{B}}^T & 0
 \end{array}\right)\,\quad\,,\quad\,
{\cal I}\,:=\, 
\left(\begin{array}{cc}\,I &0\\
 0 &\,I
 \end{array}\right)\, .
 \end{equation} 
 The operators ${\tilde{{\cal B}}}$ and ${\cal I}$ are elements of the Banach space
 ${\cal L}\hspace*{-.19cm}{{\cal L}}(l\hspace*{-.075cm}{l}_{2}\times l\hspace*{-.075cm}{l}_{2},l\hspace*{-.075cm}{l}_{2}\times l\hspace*{-.075cm}{l}_{2})$ of linear and bounded operators.  
 Now we have made all preparations to follow the arguments of Section
\ref{Secvelocity}. 
To show the compactness of ${\tilde{{\cal B}}}$ we use that ${\tilde{{\cal B}}}$ can be approximated by a sequence
 of finite operators $\{{\tilde{{\cal B}_{\ell}}}\}_{\ell \in {\bb N}}$, 
 $\{{\tilde{{\cal B}_{\ell}}}\}_{\ell \in {\bb N}}\,\in\,{\cal L}\hspace*{-.19cm}{{\cal L}}(l\hspace*{-.075cm}{l}_{2}\times l\hspace*{-.075cm}{l}_{2},l\hspace*{-.075cm}{l}_{2}\times l\hspace*{-.075cm}{l}_{2})$, with $\{{\tilde{{\cal B}_{\ell}}}\}_{\ell \in {\bb N}}\,\xrightarrow [{\cal L}\hspace*{-.14cm}{{\cal L}}]{} \,
 {\tilde{{\cal B}}}$, where we use $\forall \,\ell\,\in \,{\bb N}$  matrices
 ${\tilde{B}}_{\ell}$ \, defined along the lines of \eqref{finite} for ${\tilde{{\cal C}}}$. The convergence $\{{\tilde{{\cal B}_{\ell}}}\}_{\ell \in {\bb N}}\,\xrightarrow [{\cal L}\hspace*{-.14cm}{{\cal L}}]{} \, {\tilde{{\cal B}}}$ follows from \eqref{eststream}, \eqref{streamvollst},  $\lim_{j\rightarrow\infty} \mu_{j}=\infty\,$ 
 {\underline{or}} $\lim_{k\rightarrow\infty}\omega_{k}=\infty\,$.\,
 Finally we use Lemma \ref{LambdaNicht0} again in a modified formulation: $\lambda =0$ is not an eigenvalue of ${\tilde{{\cal B}}}$
 resp. of
problem (\ref{E1})-(\ref{E2}) and we set for the proof $ {\bw}={\underline{\nabla}}\,\Psi {{\times}} \beeis\,$.
\,\Qed
 %%%%%%%%%
\section{Orthogonality relations}
%%%%%%%%%%%%%%%%%%%
Firstly, we focus on the multiplicity of the eigenvalues.
\begin{notation} \label{multiplicity_ev} Let ${\tilde{{\cal C}}}\,\in\,
 {\cal L}\hspace*{-.19cm}{{\cal L}}(l\hspace*{-.075cm}{l}_{2}\times
 l\hspace*{-.075cm}{l}_{2},l\hspace*{-.075cm}{l}_{2}\times
 l\hspace*{-.075cm}{l}_{2})$ be the operator stated in Theorem
 \ref{T3.5} (equivalent to Problem A (\eqref{E1}-\eqref{E2}) for  
all $\mathcal{A}$). 
We call $N({\tilde{{\cal C}}} -\lambda\,{\cal I})$  as a subspace of 
$l\hspace*{-.075cm}{l}_{2}\times l\hspace*{-.075cm}{l}_{2}$ the null space of the operator ${\tilde{{\cal C}}} -\lambda\,{\cal I}$ for 
$\lambda\,\in\,{\bb C}$:
\begin{align}\label{Nullraum} 
N({\tilde{{\cal C}}} -\lambda\,{\cal I})\,:=\,\{{\tilde{\mathbf{x}}}\,\in\,
l\hspace*{-.075cm}{l}_{2}\times l\hspace*{-.075cm}{l}_{2}:\,(
{\tilde{{\cal C}}} -\lambda\,{\cal I}){\widetilde{\mathbf{x}}}\,=\,{\widetilde{\mathbf{o}}}\}\,,
\end{align}
where ${\widetilde{\mathbf{x}}}\,=\,
({\tilde{\mathbf{c}}}
 ,{\tilde{\mathbf{d}}})\,,{\widetilde{\mathbf{o}}}\,\in\,l\hspace*{-.075cm}{l}_{2}\times l\hspace*{-.075cm}{l}_{2}$, with
${\widetilde{\mathbf{o}}}$ as the zero element of the space
$l\hspace*{-.075cm}{l}_{2}\times l\hspace*{-.075cm}{l}_{2}$\,. 
The multiplicity of the eigenvalue $\lambda$ of ${\tilde{{\cal C}}}$ is
\begin{align}
\label{EV-Mult} 
\mathfrak{N}(\lambda)\,:=\,\text{dim}\,N({\tilde{{\cal C}}} -\lambda\,{\cal I})\,.
\end{align}
\end{notation}
\noindent We will prove the following orthogonality relations 
regarding the eigen-triples of Problem A:
\begin{lem}
%The eigenvalues $\lambda$ are different from zero
For any eigenvalue $\lambda$ of ${\tilde{{\cal C}}}$ (resp. of the
Problem A (cf. \eqref{E1}-\eqref{E2}), there also exists
the eigenvalue $(-\lambda)$. We see that for every eigen-triple 
$(\lambda,{\bw},{\theta})\in  
{\bb C}\times{\underline{\bb V}}\times {\W}
$ there exists a corresponding eigen-triple 
$(-\lambda,{\bw}^{-},{\theta}^{-})\,:=\,(-\lambda,-{\bw},{\theta})\,\in  
{\bb C}\times{\underline{\bb V}}\times {\W}
$. The eigenvalues $\lambda$ and $-\lambda$ have the same
multiplicity: $ \mathfrak{N}(\lambda)\,=\, 
\mathfrak{N}(-\lambda)$.
\end{lem}
\noindent {\bf Proof.} 
Let us assume that 
$(\lambda,{\bw},{\theta})\,\in\,
{\bb C}\times {\underline{\bb V}}^{2}\times {\bb
  W}_{2}^{1}\hspace{-.62cm}{~}^{{~}^{{~}^{o}}}\hspace{.2cm}
\cap {\bb W}_{2}^{2}$  
is a fixed solution of eqs.~(\ref{E1}), (\ref{E2}). Then by inserting
$(-\lambda,-{\bw},{\theta})$ in
\eqref{E1} and \eqref{E2} (resp. multiplying the second equation
\eqref{E2} with -1) we find 
\begin{eqnarray*}
\theta{\underline{\nabla}} S -\lambda\Delta(-\bw)&= &-{\underline{\nabla}} \,p\,, \quad {\underline{\nabla}}^T
                                                                            \cdot\widetilde{\bw}=0\,,\\
-\bw^T\cdot{\underline{\nabla}} S-\lambda\Delta\theta&=&0\,.
\end{eqnarray*}
It is hence evident that the problem is fulfilled  by
\[(-\bw,\theta )\,\,,\]
where the elements $(\bw,\theta )$ and $(-\bw,\theta )$ are linearly 
independent in the Cartesian product.\\
Now we use the linear and bounded operators ${\tilde{{\cal C}}}$ and ${\cal I}$ ;\quad 
$
{\tilde{{\cal C}}},\,{\cal I}: l\hspace*{-.075cm}{l}_{2}\times l\hspace*{-.075cm}{l}_{2}\longrightarrow  l\hspace*{-.075cm}{l}_{2}\times l\hspace*{-.075cm}{l}_{2}$ \eqref{OPDEF}
and the equivalent eigenvalue problem (cf. \eqref{EVPlambda}). The sequences of Fourier coefficients  ${{\mathbf{c}}}$ and ${{\mathbf{}d}}$ corresponding to ${\bw}$ and
${\theta}$ belong to the spaces 
$l\hspace*{-.075cm}{l}^{{\bb V}}_{2}$ and $l\hspace*{-.075cm}{l}^{{\bb W}}_{2}$.
So we apply the notations ${\tilde{\mathbf{c}}}$ and 
${\tilde{\mathbf{d}}}$, (${\tilde{\mathbf{c}}},\,{\tilde{\mathbf{d}}}\,\in\,l\hspace*{-.075cm}{l}_{2}$) cf. \eqref{AbbrsequWeight} for the correspondent
eigen-triple $(\lambda,{\tilde{\mathbf{c}}},{\tilde{\mathbf{d}}})\,\cong
(\lambda,{\bw},{\theta})$. 
Let us suppose that the eigen-triple $(\lambda,{\tilde{\mathbf{c}}},{\tilde{\mathbf{d}}})\,\cong
(\lambda,{\bw},{\theta})$ fulfils the equation:
 \begin{equation}\label{diskussEVPlambda}
{\tilde{{\cal C}}}\left(\begin{array}{c} {\tilde{\mathbf{c}}}\\
 {\tilde{\mathbf{d}}}
 \end{array}\right)\,=\,
\left(\begin{array}{cc} 0&{\tilde{C}}\\
 {\tilde{C}}^T & 0
 \end{array}\right)
\left(\begin{array}{r} {\tilde{\mathbf{c}}}\\
 \,{\tilde{\mathbf{d}}}
 \end{array}\right)\,=\, 
\left(\begin{array}{r} {\tilde{C}}{\tilde{\mathbf{d}}}\\
 \, {\tilde{C}}^T {\tilde{\mathbf{c}}}
 \end{array}\right)\,=\,
\lambda\,{\cal I}\left(\begin{array}{c} {\tilde{\mathbf{c}}}\\
 {\tilde{\mathbf{d}}}
 \end{array}\right)\,=\,\lambda\, \left(\begin{array}{c} {\tilde{\mathbf{c}}}\\
 {\tilde{\mathbf{d}}}
 \end{array}\right)\,.
 \end{equation}
It is obvious, that also $(-\lambda,-{\tilde{\mathbf{c}}},{\tilde{\mathbf{d}}})\,\cong
(-\lambda,-{\bw},{\theta})$ satisfies the ``ruling" \,equation as an eigen-triple.
 \begin{equation}\label{PlusMinusGLS3}
{\tilde{{\cal C}}}\left(\begin{array}{r} -{\tilde{\mathbf{c}}}\\
 \,{\tilde{\mathbf{d}}}
 \end{array}\right)\,=\,
\left(\begin{array}{cc} 0&{\tilde{C}}\\
 {\tilde{C}}^T & 0
 \end{array}\right)
\left(\begin{array}{r} -{\tilde{\mathbf{c}}}\\
 \,{\tilde{\mathbf{d}}}
 \end{array}\right)\,=\, 
\left(\begin{array}{r} {\tilde{C}}{\tilde{\mathbf{d}}}\\
 \, -{\tilde{C}}^T {\tilde{\mathbf{c}}}
 \end{array}\right)\,=\,
 -\lambda\,{\cal I}
 \left(\begin{array}{r} -{\tilde{\mathbf{c}}}\\
 \,{\tilde{\mathbf{d}}}
 \end{array}\right)\,\,=-\lambda
 \left(\begin{array}{r} -{\tilde{\mathbf{c}}}\\
 \,{\tilde{\mathbf{d}}}
 \end{array}\right)\,.
 \end{equation}  
The multiplicities of $\lambda$ and $-\lambda$ are equal. To see this we use 
\eqref{diskussEVPlambda}, \eqref{PlusMinusGLS3} and $ \mathfrak{N}(\lambda)\,<\,\infty$ from Theorem \ref{T3.5}. 
Then we choose
$({\tilde{\mathbf{c}}},{\tilde{\mathbf{d}}})\,\in\,N({\tilde{{\cal C}}} -\lambda\,{\cal I})$
as a fixed element of an $l\hspace*{-.075cm}{l}_{2}\times l\hspace*{-.075cm}{l}_{2}$-orthonormal basis of $N({\tilde{{\cal C}}} -\lambda\,{\cal I})$ to get, that $(-{\tilde{\mathbf{c}}},{\tilde{\mathbf{d}}})$ is
an element of an $l\hspace*{-.075cm}{l}_{2}\times l\hspace*{-.075cm}{l}_{2}$-orthonormal basis of $N({\tilde{{\cal C}}} + \lambda\,{\cal I})$. The assertion is showed by interchanging
the roles of $\lambda$ and $-\lambda$.
 \Qed
{~}\\[.2cm]
We write down the statement of the scalar product and orthogonality in the Cartesian product
of Hilbert spaces to illustrate the use of an orthonormal basis on $l\hspace*{-.075cm}{l}_{2}\times l\hspace*{-.075cm}{l}_{2}$ (cf. Section \ref{Secvelocity}, \eqref{EWMaxC}).\\[-.5cm]
\begin{notation} \label{Cartesian_Pr} 
Let us use the abbreviations from \eqref{Abbrsequ} for elements of the
Hilbert sequence space $l\hspace*{-.075cm}{l}_{2}$.
Then the scalar product on the Cartesian product is defined via
\begin{align*}
[(\mathbf{a},\mathbf{b}),
(\mathbf{c},\mathbf{d})]_{l\hspace*{-.075cm}{l}_{2}\times l\hspace*{-.075cm}{l}_{2}}
&:= \,\frac{1}{2}\big((\mathbf{a},\mathbf{c})_{2,(l)}\,+\,(\mathbf{b},\mathbf{d})_{2,(l)}\big)
\\
\forall\,(\mathbf{a},\mathbf{b}),\,(\mathbf{c},\mathbf{d})
&\,\in\,l\hspace*{-.075cm}{l}_{2} \times l\hspace*{-.075cm}{l}_{2}
\end{align*}
The norm on $l\hspace*{-.075cm}{l}_{2} \times l\hspace*{-.075cm}{l}_{2}$ (cf. \eqref{EWMaxC})
is taken in the usual way as
\begin{align*}
\|({\mathbf{c}},
{\mathbf{d}})\|_{l\hspace*{-.075cm}{l}_{2}\times l\hspace*{-.075cm}{l}_{2}}\,=\,\sqrt{[(\mathbf{c},\mathbf{d}),
(\mathbf{c},\mathbf{d})]_{l\hspace*{-.075cm}{l}_{2}\times l\hspace*{-.075cm}{l}_{2}}}|\,.
\end{align*}
\end{notation}
\begin{lem} \label{KartONS}
The inner product $[(\mathbf{a},\mathbf{b}),
(\mathbf{c},\mathbf{d})]_{l\hspace*{-.075cm}{l}_{2}
\times l\hspace*{-.075cm}{l}_{2}}$
given in Notation \ref{Cartesian_Pr} is a scalar product on $l_{2} \times l_{2}$ and
the Hilbert space $l\hspace*{-.075cm}{l}_{2} \times
l\hspace*{-.075cm}{l}_{2}$ is defined by 
$l\hspace*{-.075cm}{l}_{2} \times l\hspace*{-.075cm}{l}_{2}\,:=\,
(l_{2}  \times l_{2}, [(.,.),
(.,.)]_{l\hspace*{-.075cm}{l}_{2}\times l\hspace*{-.075cm}{l}_{2}})$.
Let the orthonormal basis of the $l\hspace*{-.075cm}{l}_{2}$ be 
\begin{equation}\label{deltaBasis}
\{{\text{\boldmath$\delta$}}_j\}_{j=1}^{\infty}\,\quad\text{with}\,\,
{\text{\boldmath$\delta$}}_j\,:=\,
\{{\delta}_{j,k}\}_{k=1}^{\infty}\quad \,\forall \,j\,\in\,{\bb N}\,,
\end{equation}
where ${\delta}_{j,k}$ denote the Kronecker ${\delta}$ (as function of 
the variables $j$ and $k$).\\
In this way one obtains an orthonormal basis of  $l\hspace*{-.075cm}{l}_{2}
\times l\hspace*{-.075cm}{l}_{2}$ 
by choosing the system:
\begin{equation}\label{KartBasis}
\{{\text{\boldmath$\gamma$}}_j\,=\,({\text{\boldmath$\delta$}}_j,{\text{\boldmath$\delta$}}_j)\}_{j=1}^{\infty}\quad \cup \quad \{{\text{\boldmath$\eta$}}_j\,=\,(-{\text{\boldmath$\delta$}}_j,{\text{\boldmath$\delta$}}_j)\}_{j=1}^{\infty}\,.
\end{equation}
\end{lem}
\noindent {\bf Proof.} The properties of a scalar product on $l\hspace*{-.075cm}{l}_{2} \times l\hspace*{-.075cm}{l}_{2}$
are easily checked. The only
critical point in showing that we have
an orthonormal basis of  $l\hspace*{-.075cm}{l}_{2} \times l\hspace*{-.075cm}{l}_{2}$ 
by the choosing system \eqref{KartBasis} is the simple calculation:
\begin{align*}
[{\text{\boldmath$\gamma$}}_j,{\text{\boldmath$\eta$}}_j]_{l\hspace*{-.075cm}{l}_{2}\times l\hspace*{-.075cm}{l}_{2}}\,=\,
[({\text{\boldmath$\delta$}}_j,{\text{\boldmath$\delta$}}_j),
(-{\text{\boldmath$\delta$}}_j,{\text{\boldmath$\delta$}}_j)]_{l\hspace*{-.075cm}{l}_{2}\times l\hspace*{-.075cm}{l}_{2}}\quad
\,=\,0 \,\quad\,\forall \,j\,\in\,{\bb N}\,.
\end{align*}
The completeness follows, since the zero element of $l\hspace*{-.075cm}{l}_{2} \times l\hspace*{-.075cm}{l}_{2}$ is the only element
${\text{\boldmath$\xi$}}\,\in\,l\hspace*{-.075cm}{l}_{2} \times l\hspace*{-.075cm}{l}_{2}$, 
which is orthogonal to the system \eqref{KartBasis}.
\Qed

\noindent {\bf Remark:} Important properties of the eigenvalues $\lambda$ 
of ${\tilde{{\cal C}}}$ are written down in Theorem \ref{T3.5}:
The eigenvalues $\lambda$ of ${\tilde{{\cal C}}}$ are real and $\lambda\,=\,0$ is not an
eigenvalue of  ${\tilde{{\cal C}}}$
but the only accumulation point of the eigenvalues. The set of all eigenvalues $\{\lambda\}$ 
consists of a countably infinite number of (finite multiplicity: $ \mathfrak{N}(\lambda)\,<\,\infty$) eigenvalues. The features of the construction of an
orthonormal basis in $l\hspace*{-.075cm}{l}_{2} \times l\hspace*{-.075cm}{l}_{2}$ can
be transferred to the spaces ${\underline{\bb H}}\times {\bb L}_{2}$ and 
${\underline{\bb V}}\times {\W}$ in a one to one manner. The connections of the factors
in the Cartesian products are
conserved for the Cartesian products together with the features:
${\underline{\bb V}}\,:=\,{\boldsymbol S}^{-\frac{1}{2}}({\underline{\bb H}})$ and 
${\W}\,:=\,{\boldsymbol L}^{-\frac{1}{2}}( {\bb L}_{2})$. There we have considered the Stokes operator and the Laplacian as positive definite operators in the context of spectral operators.
\begin{prop}\label{PropOrth} ({\em{Orthogonality}}) Let $\lambda_1,\lambda_2$ \,be
eigenvalues of\,
${\tilde{{\cal C}}}\,\in\,
 {\cal L}\hspace*{-.19cm}{{\cal L}}(l\hspace*{-.075cm}{l}_{2}\times l\hspace*{-.075cm}{l}_{2},l\hspace*{-.075cm}{l}_{2}\times l\hspace*{-.075cm}{l}_{2})$,
with $\lambda_1\,\not=\,\lambda_2$, then the corresponding eigenspaces
$N({\tilde{{\cal C}}} -\lambda_1\,{\cal I})$ and $N({\tilde{{\cal C}}}
-\lambda_2\,{\cal I})$ are 
orthogonal subspaces in $l\hspace*{-.075cm}{l}_{2}\times l\hspace*{-.075cm}{l}_{2}$.
\end{prop}
\noindent{\bf Proof.} According to Theorem \ref{T3.5} ${\tilde{{\cal
      C}}}$ is a self-adjoint 
operator on $l\hspace*{-.075cm}{l}_{2}\times l\hspace*{-.075cm}{l}_{2}$ and $\lambda
\,=\,0$ is not an eigenvalue of ${\tilde{{\cal C}}}$.
The statement is obtained by the application of 
 \cite[Sec. 4.2.3]{Triebel}.
\Qed
\begin{theor} \label{CompleteI}{\em{(Completeness)}} 
The set of orthonormalized eigenfunctions of ${\tilde{{\cal C}}}$ generates a complete basis
for the Cartesian product
space $l\hspace*{-.075cm}{l}_{2}\times l\hspace*{-.075cm}{l}_{2}$. If by 
$\{\lambda_{\ell}\}_{{\ell}=1}^{\infty}$ we denote the set of eigenvalues of ${\tilde{{\cal C}}}$,
we see
\begin{align}\label{NulluCart}
l\hspace*{-.075cm}{l}_{2}\times l\hspace*{-.075cm}{l}_{2}\,\,=\,{\overline{
{\bigcup_{{\ell}=1}^{\infty}}N({\tilde{{\cal C}}} -\lambda_{\ell}\,{\cal I})}}^{l\hspace*{-.075cm}{l}_{2}\times l\hspace*{-.075cm}{l}_{2}}\,\,.
\end{align}
\end{theor}
\noindent{\bf Proof.} For every $\lambda_{\ell}$ there exists a finite orthonormal basis in
$N({\tilde{{\cal C}}} -\lambda_{\ell}\,{\cal I})$.  According to Proposition
\ref{PropOrth} the eigenspaces for different $\lambda$ are orthogonal
to each other.
So we have only to show \eqref{NulluCart}. We use, that $\lambda\,=\,0$ is not an eigenvalue of 
${\tilde{{\cal C}}}$ and that the operator ${\tilde{{\cal C}}}$
is self-adjoint and compact. So we apply the orthogonal decomposition Theorem 
and the Theorem about the spectrum of compact operators from 
\cite[Sec. 4.5.2]{Triebel}) to get the justification of \eqref{NulluCart}.
\Qed

\begin{prop} \label{propLAL} 
Let $\LAL \colon D({\boldsymbol S})\times D({\boldsymbol L}) \to {\underline{\bb H}}\times {\bb L}_{2}$ be regarded as in \eqref{LAL}. Then is $\LAL$ a self-adjoint positive definite operator.
In the 
sense of \eqref{LAL} and of the matrix representation on product
spaces $\LAL$ is a component-by-component combination of the
self-adjoint positive definite operators ${\boldsymbol S}$ and
${\boldsymbol L}$ with 
pure point spectrum. We get for the ``diagonal matrices" 
\begin{align}\label{PotenzenLAL}
{\LAL}^{-1}\,=\,
\left(
 \begin{array}{cc}
 {\boldsymbol S}^{-1} & \underline{0} \\
 & \\
 ({\underline{0}})^T \,
& {\boldsymbol L}^{-1}
 \end{array}
 \right)
 \,,\,\,
 {\LAL}^{-\frac{1}{2}}\,=\,
\left(
 \begin{array}{cc}
 {\boldsymbol S}^{-\frac{1}{2}} & \underline{0} \\
 & \\
 ({\underline{0}})^T \,
& {\boldsymbol L}^{-\frac{1}{2}}
 \end{array}
 \right)\,\,\,\text{and}\,\quad {\LAL}^{\frac{1}{2}}\,=\,
\left(
 \begin{array}{cc}
 {\boldsymbol S}^{\frac{1}{2}} & \underline{0} \\
 & \\
 ({\underline{0}})^T \,
& {\boldsymbol L}^{\frac{1}{2}}
 \end{array}
 \right)\,,
\end{align}
and for
${\underline{\mathfrak{u}}}:=(\underline{\Psi},\varphi)\,\in\,{\underline{\bb
    V}}\times \W$ 
\begin{align}\label{TrickLAL}
\left({\LAL}^{\frac{1}{2}}{\underline{\mathfrak{u}}},{\LAL}^{\frac{1}{2}}{\underline{\mathfrak{w}}}\right)
\,=\, ({\boldsymbol S}^{\frac{1}{2}}\underline{\Psi}, {\boldsymbol S}^{\frac{1}{2}} \bw)
\,+\, ({\boldsymbol L}^{\frac{1}{2}}\varphi, {\boldsymbol L}^{\frac{1}{2}} \theta)\,=\,
\left(\underline{\Psi},\bw\right)_{D}\,+\,
\left(\varphi,\theta\right)_{D}\,.
\end{align}
The energetic space of \,${\LAL}$ \, is ${\underline{\bb V}}\times \W\,=\,
{\LAL}^{-\frac{1}{2}}({\underline{\bb H}}\times {\bb L}_{2})$. \,
The self-adjoint operators \,${\LAL}^{-1}$ and \,${\LAL}^{-\frac{1}{2}}$ are compact in  
${\cal L}\hspace*{-.19cm}{{\cal L}}({\underline{\bb H}}\times {\bb L}_{2},{\underline{\bb H}}\times {\bb L}_{2})$. This is also true for \,${\LAL}^{-1}$ regarded as:\,
${\LAL}^{-1}\colon {\underline{\bb H}}\times {\bb L}_{2}\,\to \,
{\underline{\bb V}}\times \W \,$.
\end{prop}
\noindent {\bf Proof} is a one to one transference of 
\cite[Sec. 4.5.2]{Triebel}  on product spaces.\Qed
\begin{theor} \label{L2CompTh} For all ${\cal A}$ Problem A 
is equivalent to the eigenvalue problem of the
compact self-adjoint operator $\LCL$ on the Cartesian product
space ${\underline{\bb H}}\,\times
{\bb L}_{2}$.
\end{theor}
\noindent {\bf Proof.} Let us use Notation \ref{SymmOp} and \eqref{SymmDefi}. So we write
\begin{equation}\label{VollstII}
 \LL ((\bw,\theta))\,=\,(\LSL \,- \,{\lambda}\LAL) ((\bw,\theta))\,=\,
 \left(\left(
\begin{array}{cc}
0  & {\underline{\nabla}}\;\!S\,\\
 &\\
({\underline{\nabla}}\;\!S)^T \,
& 0
\end{array}\right)
\,-\,{\lambda}
\left(
 \begin{array}{cc}
 {\boldsymbol S} & \underline{0} \\
 & \\
 ({\underline{0}})^T \,
& {\boldsymbol L}
 \end{array}
 \right)
 \right)
 \left(\begin{array}{c}
 \bw\\  \\  \theta
 \end{array}
 \right) ,
 \end{equation}
 where one can understand the zeros in $\LSL$ also as the zero elements
 of  ${\cal L}\hspace*{-.19cm}{{\cal L}}({\underline{\bb H}},{\underline{\bb H}})$
  resp. 
 ${\cal L}\hspace*{-.19cm}{{\cal L}}({\bb L}_{2},{\bb L}_{2})$.\\
The operators $\LAL \colon D({\boldsymbol S})\times D({\boldsymbol L}) \to {\underline{\bb H}}\times {\bb L}_{2}$ \, and
$\LSL\colon {\underline{\bb H}} \times {\bb L}_{2} \to 
{\underline{\bb H}} \times {\bb L}_{2}$ in \eqref{VollstII} are also defined through the
%via $\glqq$${\bb L}_{2}$\grqq \ \\
bilinear forms at $(\underline{\Psi},\varphi)\,\in\,
{\underline{\bb V}} \times \W $:
\begin{align}\label{LAL}
((\underline{\Psi},\varphi),\LAL(\bw,\theta))&:=\,\left(\underline{\Psi},\bw\right)_{D}\,+
\left(\varphi,\theta\right)_{D}
\\
((\underline{\Psi},\varphi),\LSL(\bw,\theta))&:=\,\left({\underline{\Psi}},\theta\,
{\underline{\nabla}}\;\!S \right)\,+\,
\left(\varphi,({\underline{\nabla}}\;\!S)^T \bw\right)
\nonumber.
\end{align}
Eqs. \eqref{LAL} transfer in the sense on dense definition. 
First we investigate the properties of $\LAL$ to establish a relationship to the 
task of Problem A.
The nontrivial solutions $(\lambda,{\bw},\theta)\in {\bb C}\times {\underline{\bb V}}\times \W$
of (\ref{E1})-(\ref{E2}) are equations in ${\underline{\bb V}}'$ 
resp. ${\bb W}_{2}^{-1}$. 
Using the notation ${\LAL}_{en}$ for the energetic expansion of ${\LAL}$ in the 
sense of  ${\LAL}_{en} \colon {\underline{\bb V}}\times \W\,\to \,{\underline{\bb V}}'
\times {\bb W}_{2}^{-1}$ and
${\underline{\mathfrak{w}}}:=(\bw,\theta)\,\in\,{\underline{\bb V}}\times \W$ we may write%##
\begin{align}\label{numberA1}
\left({\LSL}
\,-\lambda {\LAL}_{en}\right){\underline{\mathfrak{w}}}\,=\,{\underline{0}}\,\,.
\end{align} 
We prove the boundedness of ${\LSL}$
using $\gamma_{S}$ (cf. \eqref{symm1}). The crucial tool for our proof 
is the mapping property: By an appropriate norm on ${\underline{\bb V}}\times \W $
we have ${\LAL}^{\frac{1}{2}}$ as a unitary map on 
${\underline{\bb H}}\times {\bb L}_{2}$ etc. The idea of one to one mappings result in:
${\underline{\mathfrak{u}}}\,=\,{\LAL}^{-\frac{1}{2}}{\underline{\mathfrak{v}}}$ and
${\underline{\mathfrak{w}}}\,=\,{\LAL}^{-\frac{1}{2}}{\underline{\mathfrak{y}}}$. 
For Problem A we find using \eqref{VollstII} and
\eqref{LAL} that
\begin{align}\label{LCL}
({\LAL}^{-\frac{1}{2}}{\underline{\mathfrak{v}}},\LSL \,{\LAL}^{-\frac{1}{2}}{\underline{\mathfrak{y}}})\,-\, \lambda ({\LAL}^{\frac{1}{2}}{\LAL}^{-\frac{1}{2}}{\underline{\mathfrak{v}}},{\LAL}^{\frac{1}{2}}{\LAL}^{-\frac{1}{2}}{\underline{\mathfrak{y}}})\,
& =\,\nonumber\\
({\underline{\mathfrak{v}}},{\underbrace{
{\LAL}^{-\frac{1}{2}}\LSL \,{\LAL}^{-\frac{1}{2}}}_{=: \LCL}}
{\underline{\mathfrak{y}}})\,-\, \lambda ({\underline{\mathfrak{v}}},{\underline{\mathfrak{y}}})
\,\,& =\,\,0\,.
\end{align}
Finally we use the identity $\LIL$ in ${\underline{\bb H}}\times {\bb L}_{2}$ to get
\begin{align}\label{CFERTIG}
(\LCL \,-\,\lambda \LIL) {\underline{\mathfrak{y}}}\,=\,{\underline{0}}
\end{align}
as the equivalent eigenvalue problem in ${\underline{\bb H}}\times {\bb L}_{2}$
for the self-adjoint compact operator $\LCL$.
\Qed
%\newpage
% {~}\\[.4cm]}}Remainder:\\
% {\red{One has also as a consequence of the former considerations
% to orthogonality relations using fixed eigen-triples $(\lambda,{\bw},{\theta})$ 
% and $(-\lambda,-{\bw},{\theta})$ with $({\bw},{\theta})\,\cong\,({\tilde{\mathbf{c}}},{\tilde{\mathbf{d}}})\,\in\,N({\tilde{{\cal C}}} -\lambda\,{\cal I})$
% resp. $(-{\bw},{\theta})\,\cong\,(-{\tilde{\mathbf{c}}},{\tilde{\mathbf{d}}})\
% \in\,N({\tilde{{\cal C}}} +\lambda\,{\cal I})$, 
% \begin{displaymath}
% \begin{array}{l|l|}
% \cline{2-2}{~}&  {~}\\
% {~}&
% \|\bw\|_{D}^2\,=\, \|\theta\|_{D}^2\,> \,0 \quad \quad \text{and}\\ 
% \text{that}
% \qquad& \\
% {~}& \lambda \cdot (\|\bw\|_{D}^2\,+\, \|\theta\|_{D}^2)\,=\,2 \cdot (\theta \:\!
% {\underline{\nabla}}\:\! S ,\bw)\,.
%  \\ 
% \cline{2-2}
% \end{array}
% \end{displaymath}}}
%\color{black}

%\input{append.tex}

\end{document}